\documentclass{article}
\usepackage{amssymb}
\usepackage{amsmath}
\usepackage{graphicx}
\usepackage{color}
\usepackage{subfloat}
\usepackage{subcaption}
\usepackage{cite}
\usepackage{todonotes}

\usepackage{hyperref}			

\usepackage{amsthm}               
\usepackage{amsmath}             

\newtheorem{teo}{Theorem}[section]
\newtheorem{defi}[teo]{Definition}
\newtheorem{lem}[teo]{Lemma}

\newtheorem{ass}[teo]{Assumption}

\newcommand{\REV}[1]{#1}

\newcommand{\metls}{{\tt SGN\_RC }}
\newcommand{\metlspunto}{{\tt SGN\_RC}}

\newcommand{\metsist}{{\tt SGN\_JS }}
\newcommand{\metsistp}{{\tt SGN\_JS}}
\newcommand{\be}{\begin{equation}}
\newcommand{\ee}{\end{equation}}
\newcommand{\ba}{\begin{array}}
\newcommand{\ea}{\end{array}}
\newcommand{\bee}{\begin{eqnarray*}}
\newcommand{\eee}{\end{eqnarray*}}
\newcommand{\bea}{\begin{eqnarray}}
\newcommand{\eea}{\end{eqnarray}}

\newcommand{\R}{\mathbb{R}}
\def\sc{\mathsf{oc}}
\def\toc{\mathsf{TOC}}

\newcommand{\comment}[1]{}

\newcommand{\E}{\mathbb E}
\newcounter{algo}[section]
\renewcommand{\thealgo}{\thesection.\arabic{algo}}
\newcommand{\algo}[3]{\refstepcounter{algo}
\begin{center}
\framebox[\textwidth]{
\parbox{0.95\textwidth} {\vspace{\topsep}
{\bf Algorithm \thealgo . #2}\label{#1}\\
\vspace*{-\topsep} \mbox{ }\\
{#3} \vspace{\topsep} }}
\end{center}}
\def\IR{\hbox{\rm I\kern-.2em\hbox{\rm R}}}
\def\IC{\hbox{\rm C\kern-.58em{\raise.53ex\hbox{$\scriptscriptstyle|$}}
    \kern-.55em{\raise.53ex\hbox{$\scriptscriptstyle|$}} }}
    
\newcommand{\RR}{\mathbb{R}}

\def \E {{\mathbb E}}

\usepackage{bbm}				
\def \1{{\mathbbm 1}}

\newcommand{\mqp}{{\tt NLS}}	    
\newcommand{\sistemap}{{\tt NSE}}
\newcommand{\mq}{{\tt NLS }}	    
\newcommand{\sistema}{{\tt NSE }}

\newcommand{\cS}{{\cal S}}

\DeclareMathOperator{\spa}{span}
\DeclareMathOperator{\diag}{diag}
\DeclareMathOperator{\argmin}{argmin}

\newcommand{\minor}[1]{#1}
\newcommand{\tblue}[1]{#1}

\title{Inexact Gauss-Newton methods with matrix approximation by sampling for nonlinear least-squares and systems
}
 \author{ Stefania Bellavia\footnotemark[1], Greta Malaspina\footnotemark[1], Benedetta Morini\footnotemark[1]}

\date{}

\begin {document}
\maketitle
\footnotetext[1]{Dipartimento  di Ingegneria Industriale, Universit\`a degli Studi di Firenze,
Viale G.B. Morgagni 40,  50134 Firenze,  Italia. Members of the INdAM Research Group GNCS. Emails:
stefania.bellavia@unifi.it, greta.malaspina@unifi.it,  benedetta.morini@unifi.it}
\footnotetext[2]{ The research that led to the present paper was partially supported by INDAM-GNCS through Progetti di Ricerca 2023 and 
by PNRR - Missione 4 Istruzione e Ricerca - Componente C2 Investimento 1.1, Fondo per il Programma Nazionale di Ricerca e Progetti di Rilevante Interesse Nazionale (PRIN) funded by the European Commission under the NextGeneration EU programme, project ``Advanced optimization METhods for automated central veIn Sign detection in multiple sclerosis from magneTic resonAnce imaging (AMETISTA)'',  code: P2022J9SNP,
MUR D.D. financing decree n. 1379 of 1st September 2023 (CUP E53D23017980001), project 
``Numerical Optimization with Adaptive Accuracy and Applications to Machine Learning'',  code: 2022N3ZNAX
 MUR D.D. financing decree n. 973 of 30th June 2023 (CUP B53D23012670006).}

\begin{abstract}
We develop and analyze stochastic inexact Gauss-Newton methods for nonlinear least-squares problems and for  nonlinear systems of
equations.
Random models are formed using suitable sampling strategies for the matrices involved in the deterministic models. The analysis of the expected number of iterations needed in the worst case to achieve a desired level of accuracy in the first-order optimality condition provides guidelines for applying sampling and enforcing, with \minor{a} fixed probability, a suitable accuracy in the random approximations.
Results of the numerical validation of the algorithms are presented.

\end{abstract}

\section{Introduction}
This work addresses the solution of large-scale nonlinear least-squares problems and nonlinear systems by  inexact Newton  methods \cite{DES} combined with  random models and the line-search strategy.
The Nonlinear Least-Squares problem (\mqp) has the form
\begin{equation}\label{ls}
\min_{x\in\RR^n} f(x)= \frac{1}{2m} \|R(x)\|_2^2,
\end{equation}
with $R:\RR^n\rightarrow \RR^m$, $m\ge n$, continuously differentiable. 
As a special case, problem (\ref{ls}) includes  the solution of the square Nonlinear System of Equations (\sistemap)
\be \label{sistema}
F(x)=0,
\ee
where $F\colon \RR^n\rightarrow\RR^n$ is continuously differentiable; in fact,
the solutions of the nonlinear system are zero-residual solutions of the problem 
\be\label{problem1}
\min_{x\in\RR^n} f(x) = \frac{1}{2}\| F(x) \|_2^2.
\ee
In case $F=\nabla G$ and $G$ is an {\em invex} differentiable function, solving (\ref{sistema}) is equivalent to minimizing $G$ \cite{invex}.

In the following, we will refer to functions $R$ and $F$ as residual functions, irrespective of the problem under consideration, and the form (\ref{ls}) 
or (\ref{problem1}) of $f$ will be understood from the context.

To improve the computational complexity of the  deterministic inexact Newton methods, the procedures presented here use  models inspired by randomized linear algebra, see e.g., \cite{acta, tropp}\minor{.}  Indeed, random  approximations of expensive  derivatives, such as the gradient of $f$ and the Jacobian of the residual function, and random approximations of the Jacobian-vector product  can considerably reduce the computational effort of the solvers \cite{bg, BGM,bkk, bkm, sirtr_coap, ferrara, bbn_2017,bbn,STORM1,  acta,mit1,xu2}. \minor{We address this issue by using matrix approximation through sampling.}
Specifically, at each iteration, three major tasks are performed. The first task consists in building a random linearized model of the residual function; sampling is used for approximating the Jacobian of the residual function and our approach includes random compression, random  sparsification and standard batch approximations; the gradient of $f$ is easily obtained  as a byproduct.
The resulting procedures  fall in the class of stochastic Gauss-Newton type methods. The second task is the  approximate minimization of the model via a proper Krylov solver in order to compute the inexact step. The third task is the test acceptance of the trial step by means of the Armijo condition; function $f$ is supposed to be evaluated exactly while the gradient of $f$ is random. 
We discuss the strategies for  building the random models, outline the accuracy requests made with some probability for such models and obtain a bound on the expected number of performed iterations to achieve a desired level of accuracy in the first-order optimality condition. Further, we provide a preliminary numerical validation of our algorithms.

The recent literature on optimization with random or noisy models  is vast, restricting to line-search approaches  some recent contributions are
\cite{ferrara, Berahas, bbn, cartis, PS}.  Referring to problems (\ref{ls}) and (\ref{problem1}) and line-search and/or Inexact Newton  methods we are aware of  papers \cite{LN, LR, Newt-MR1, sketched, wwz2023, WCK}. 
In particular, the  deterministic method for (\ref{problem1})  proposed in \cite{LN} is based on a sparsification of the Jacobian. 
The inexact Newton-Minimal residual methods  proposed in \cite{LR,Newt-MR1} employ  exact function evaluations and
are applicable to problem (\ref{sistema}) if the Jacobian is symmetric; the exact Jacobian matrix is used in  \cite{Newt-MR1}
while  approximations of the  Jacobian under a  deterministic and uniform accuracy requirement are used in \cite{LR}.  Globalization strategies are not discussed in   \cite{LR, Newt-MR1}.
Paper \cite{sketched} relies on sketching matrices to reduce the dimension of the Newton system for possibly non-square nonlinear systems. The method in \cite{wwz2023} is a stochastic regularized Newton method for (\ref{sistema})  with batch approximations for the function and the Jacobian, the resulting trial step is used if it can be accepted by an inexact line-search condition, otherwise a preset step is taken.
Finally, the method in \cite{WCK} is a locally convergent Newton-GMRES method for Monte-Carlo based mappings.
Comparing with  \cite{LN, LR, Newt-MR1, sketched}, our contribution consists of the use of approximations of the derivatives based on  randomized linear algebra, and of globalization via line-search, and includes deterministic and adaptive accuracy in the limit case where accuracy requirement are  satisfied almost surely.
Comparing with \cite{wwz2023, WCK}, we address approximation of the Jacobian matrix via sparsification besides the considered mini-batch approximation.

The paper is organized as follows. In Section 2 we introduce and discuss our algorithms for \mq and \sistema problems, in Section 3 we perform the theoretical analysis and obtain the expected number of
iterations required to reach an approximate first-order optimality point, in Section 4 we 
present  preliminary  numerical results for our algorithms and in Section 5 we give some conclusions.

\subsection{Notations}
We denote the $2$-norm as $\|\cdot\|$,  the infinity norm as $\|\cdot\|_{\infty}$, the Frobenius norm as $\|\cdot\|_F$. The components of the residual functions
are denoted as $R(x)=\left( R(x)_{(1)},  \cdots, R(x)_{(m)} \right )^T$, $F(x)=\left( F(x)_{(1)},  \cdots, F(x)_{(n)} \right )^T$.
The Jacobian matrix of both $R$ and $F$ is denoted as $J$ with dimension specified by the problem, i.e., 
$J\colon \RR^{n}\rightarrow\RR^{m\times n}$ for \mq problem and
$J\colon \RR^{n}\rightarrow\RR^{n\times n}$ for \sistema problem. 
The probability of an event is denoted as $\mathbb P({\rm event})$,
and $\1({\rm event})$ is the indicator function of the event.

\section{Inexact line-search methods}
We introduce the general scheme of our procedure and then specialize the construction of the random model and the computation of the step\minor{s.}

The $k$-th iteration of our method is sketched in Algorithm \ref{algo_general}.
Given $x_k\in {\mathbb R}^n$ and the positive step-length $t_k$, we linearize  the residual function at $x_k+s$ and build a random model $\widetilde m_k(s)$ which replaces the deterministic model 
\begin{eqnarray}
m_k(s) &=&  \frac{1}{2m}\|J(x_k)s+R(x_k)\|^2, \label{model_nsl}\\
m_k(s)& = &  \frac 1 2 \|J(x_k)s+F(x_k)\|^2, \label{model_nse}
\end{eqnarray}
for \mq problem and \sistema problem, respectively.
Along with $\widetilde m_k$ we compute a stochastic approximation $g_k$ of the gradient $\nabla f(x_k)$.

The tentative step $s_k$ is then computed minimizing $\widetilde m_k$ in a suitable subspace ${\cal{K}}_{k}^{(\REV{\ell})}$ of $\RR^n$
\begin{equation}\label{stepin}
s_k=\argmin_{s\in {\cal{K}}_{k}^{(\REV{\ell})} } \widetilde m_k(s).
\end{equation}

Once $s_k$ is available,  we  test the Armijo condition  (\ref{armijo}) using exact evaluations of $f$ and the stochastic gradient $g_k$. If $x_k+t_ks_k$ satisfies such condition we say that the iteration is successful,  
accept the step and increase the step-length $t_k$ for the next iteration. Otherwise, the iteration is declared unsuccessful, 
the step is rejected and  the step-length $t_k$ is reduced for the next iteration.

 \algo{}{General scheme: $k$-th iteration }{
\noindent
Given  $ x_k\in {\mathbb R}^n, \;  c \in (0,1)$, 
$\tau \in (0,1)$,  $t_{\max}>0$, $t_k \in(0, t_{\max}]$.
\vskip 2 pt
\begin{description} 
\item{Step 1.} Form a random model $\widetilde m_k(s)$ and  the stochastic gradient $g_k$.
\\
\hspace*{15pt}
Compute $s_k$ in (\ref{stepin})   {such that $s_k^Tg_k<0.$}
\item{Step 2.} If $t_k$ satisfies condition
\be \label{armijo}
 f(x_k + t_k s_k) \leq  f(x_k) + c t_k s_k^T  g_k,
\ee
\hspace*{15pt}
Then (successful iteration)  \\
\hspace*{30pt}$ x_{k+1} = x_{k} + t_k s_k$,\quad  $t_{k+1}=\min\{t_{\max},\tau^{-1} t_k\}$,\quad $ k = k+1$
\\
\hspace*{15pt} Else (unsuccessful iteration)
\\
\hspace*{25pt } $ x_{k+1} = x_{k}$,\quad $t_{k+1}=\tau t_k$,\quad $ k = k+1$.
\end{description}
}\label{algo_general}
\vskip 5pt

In the following two sections we describe how we realize Step 1 for the
problems of interest.

\subsection{\mq problem: inexact Gauss-Newton method with row compression of the Jacobian}\label{secmodmq}
{   In this section we consider the nonlinear least-squares problem (\ref{ls})  and present an inexact procedure for building the trial step. Our approach  is based on a  random model of reduced dimension with respect to the dimension $m$ of the linear residual $(J(x_k)s+R(x_k))$ in (\ref{model_nsl}).}
A weighted random row compression is applied to $J(x_k)$; as a result $\widetilde J_k\in \mathbb{R}^{d\times n}$, $d\le m$, is formed by selecting a subset of $d$ rows of $J(x_k)$ out of $m$ and multiplying each selected row by a suitable weight. 
As for the residual function, the vector $\widetilde R_k\in \mathbb{R}^{d}$ is  formed  by selecting the  subset of $d$ entries associated to the  rows of  $\widetilde J_k$.
A practical way to form $\widetilde J_k$ and $\widetilde R_k$  is described at the end of this section.

The resulting  model
\begin{equation}\label{modelLS}
\textstyle\widetilde m_k(s)=\frac{1}{2m}\|\widetilde J_k s+ \widetilde R_k\|^2 
\end{equation}
can be approximately minimized  using an iterative method. In fact, the  $\widetilde m_k(s)$ can be \REV{minimized in nested Krylov subspaces by procedures such as LSMR \cite{LSMR} and LSQR \cite{lsqr}\footnote{See https://web.stanford.edu/group/SOL/home\_software.html}.}
These methods
can compute a minimum-length solution of $\min_s \widetilde m_k(s)$. Starting from the null initial guess $s_k^{(0)}=0$, a sequence of iterates
$\{ s_k^{(\ell)}\}$, $\ell\ge 1$, satisfying
\begin{equation}\label{lsqr}
\textstyle\|\widetilde J_k s_k^{(\ell)}+\widetilde R_k\|^2=\min_{s\in K_k^{(\ell)}} \|\widetilde J_k s+\widetilde R_k\|^2  ,
\end{equation}
is generated 
with $$K_k^{(\ell)}=\spa\left\{\widetilde J_k^T \widetilde R_k,  (\widetilde J_k^T\widetilde J_k)\widetilde J_k^T \widetilde R_k,\ldots, (\widetilde J_k^T\widetilde  J_k)^{\ell-1}   \widetilde J_k^T \widetilde R_k\right\}$$
for some integer $\ell\ge 0$. As a stopping criterion we use 
 \begin{equation}\label{stoplsqr}
 \|\widetilde J_k ^T r_k\| \le \eta_k \|\widetilde J_k^T \widetilde R_k
\|  \, \, \mbox{ with }\, \,  \textstyle r_k= \widetilde J_k s_k+ \widetilde R_k,\end{equation}
and $\eta_k\in \REV{[0, \bar \eta)}$, $\bar \eta<1$, named forcing term \cite{DES}.
We summarize this procedure in the following algorithm.
\algo{}{Step 1 of Algorithm \ref{algo_general} for \mq}{
 \noindent
Given  $ x_k\in {\mathbb R}^n$.
\vskip 2pt 
\begin{description}
\item{Step 1.1} Choose $\eta_k\in \REV{[0,\bar \eta)}$, $d\in \mathbb{N}$, $ 1\le d \le m$.\\
\hspace*{25pt}
Form $\widetilde J_k\in\RR^{d\times n}$, $\widetilde R_k\in\RR^{d}$ and $g_k=\frac{1}{m}\widetilde J_k^T\widetilde R_k$.
\item{Step 1.2}   {Apply a Krylov method to $\min_{s}\widetilde m_k(s)$ 
 with  $\widetilde m_k(s)$ given in (\ref{modelLS})
\\
\hspace*{25pt} and compute $s_k$ satisfying (\ref{stoplsqr})}.
\end{description}
}\label{algo}

\begin{lem} \label{discesa_ls}
Let $s_k, \widetilde J_k, g_k$ as in Algorithm \ref{algo}.
Then
$
s_k^Tg_k\le 0.
$
\end{lem}
\begin{proof}
By construction,  $s_k=s_k^{(\ell)}$ \REV{satisfies \eqref{lsqr} for some $\ell\ge 1$, then   the residual vector $ r_k$ in (\ref{stoplsqr}) is orthogonal to any vector in $\widetilde J_k  K_k^{(\ell)} $.} Consequently, $s_k^T \widetilde J_k^T r_k=0$ and 
\begin{equation}\label{inmq}
\frac{1}{m}s_k^T \widetilde J_k^T r_k=
\frac{1}{m}s_k^T \widetilde J_k^T \left( \widetilde J_k s_k +\widetilde R_k\right)= \REV{ \frac{1}{m} s_k^T \widetilde J_k^T  \widetilde J_k s_k+s_k^Tg_k }=0.
\end{equation}
The thesis follows since $\widetilde J_k^T \widetilde J_k$ is symmetric positive semidefinite.
\end{proof}
We conclude this section discussing the construction of the random model.
We form the matrix $\widetilde J_k$  by sampling the rows of $J(x_k)$ and the  vector $\widetilde R_k$ by sampling the components of  $R(x_k)$ accordingly.
We can build $\widetilde J_k$ and $\widetilde R_k$ as a byproduct of the gradient approximation following \cite[\S 7.3.2]{acta}.
In particular, denoting the $i$-th row of $ J(x_k)$ as $ J(x_k)_{(i,:)}$ and the $i$-th component of $R(x_k)$ as $ R(x_k)_{(i)}$, the gradient $\nabla f(x_k)$ can be expressed as
$$\nabla f(x_k) = \frac{1}{m}J(x_k)^T R(x_k) =\frac{1}{m}\sum_{i=1}^{m} (J(x_k)_{(i,:)})^T R(x_k)_{(i)}.$$
Let $p_1^k,\dots,p_m^k$ be a  probability distribution associated to   $(J(x_k)_{(i,:)})^T R(x_k)_{(i)}$, $i=1,\dots,m$,
and let ${\cal{M}}_k\subset \{1, \ldots,m\}$  be  a random subset of indices  such that index $i$ is chosen with probability $p_i^k$. We define $\overline{J}_k\in\R^{m\times n}$ as the matrix whose $i$-th row is such that
\be \label{barJ-NLS}
(\overline{J}_k)_{(i,:)} = \begin{cases} \frac{1}{|\mathcal{M}_k|p^k_i}(J(x_k)_{(i,:)}) & \text{if}\ i\in\mathcal{M}_k\\
0 & \text{otherwise}\end{cases}, 
\ee
and denote with $\widetilde{J}_k\in\R^{|\mathcal{M}_k|\times n}$ the compressed matrix obtained by retaining  the rows of $\overline{J}_k$ that correspond to indices in $\mathcal{M}_k$. 
{  Analogously, we let $\overline{R}_k\in\R^{m}$ be the vector whose $i$-th component is such that
$$(\overline{R}_k)_{(i)} = \begin{cases} R(x_k)_{(i)} & \text{if}\ i\in\mathcal{M}_k\\
0 & \text{otherwise}\end{cases}, $$
and denote $\widetilde{R}_k\in\R^{|\mathcal{M}_k|}$ the compressed vector obtained by retaining the rows of $\overline{R}_k$ that correspond to indices in $\mathcal{M}_k$.} 
We remark that $\overline{J}_k$  is an unbiased estimator of the Jacobian $J(x_k)$ and that $\widetilde{J}_k = S_kP_k^{-1}J(x_k),$ and \REV{ $\widetilde R_k=S_k R_k$} with $P_k = |\mathcal{M}_k|\diag(p^k_1,\dots,p^k_m),$  and $S_k\in\R^{|\mathcal{M}_k|\times m}$ being a suitable submatrix of the identity matrix of dimension $m$.  

A stochastic approximation of $\nabla f(x_k)$ can then be defined as
\begin{equation}\label{gk_mq} 
g_k = \frac{1}{m}\widetilde J_k^T \widetilde R_k = \frac{1}{m}\overline{J}_k^TR_k= \frac{1}{m|{\cal{M}}_k|}\sum_{i\in {\cal{M}}_k}  \frac{1}{p_i^k}(J(x_k)_{(i,:)})^T R(x_k)_{(i)}. 
\end{equation}
As for probabilities, they can be uniform, i.e., $p_i^k = 1/m$,  $i=1,\dots,m$,  or 
correspond to the so-called importance sampling \cite[\S 7.3]{acta}.
The Bernstein inequality \cite[Th. 7.2]{acta} indicates how large the cardinality of $|{\cal{M}}_k|$ should be to ensure
\begin{equation}\label{accg}
    \textstyle\mathbb P\left(\left\|\nabla f(x_k) - g_k\right\|\leq \rho\right)
\geq 1-\delta_g,
\end{equation}
given an accuracy requirement $\rho>0$ and a probability  $\delta_g\in(0,1)$.  
A general formulation of the  Bernstein inequality is given below.
\begin{teo}\cite[Th. 7.2]{acta}\label{Bernth}
Let $B\in \mathbb{R}^{q_1\times q_2}$ be a fixed matrix and let the random matrix $X\in \mathbb{R}^{q_1\times q_2}$ satisfy $\E[X]=B$ and $\|X\|\le M_X$.
Define the per-sample second moment
$v(X)=\max\left \{ \|\E\left [X^TX\right]\|, 
\|\E\left [X X^T\right ] \|\right \}$.
Form the matrix sampling estimator 
$\overline{X}_w=\frac 1 w \sum_{i=1}^{w} X_i $, where {  $w$ is an integer and} $X_i$ are 
i.i.d and have the same distribution as $X$.
Then, for all $\rho>0$
$$
   \textstyle\mathbb P\left(\left\|B-\overline{ X}_w\right\|\leq \rho\right)
\geq 1-\delta,
$$
if 
$$
w\ge \left( \frac{2v(X)}{\rho^2}+\frac{4M_X}{3\rho}\right) \log\left(\frac{q_1+q_2}{\delta}\right) .
$$
\end{teo}

Summarizing, the cardinality of the set ${\cal{M}}_k$
can be ruled by the accuracy requirement in probability specified above; once the set ${\cal{M}}_k$ is chosen, $\widetilde J_k\in \RR^{|{\cal{M}}_k|\times n}$ consists of the rows of $J(x_k)$ with index $i\in {\cal{M}}_k$,  multiplied by  suitable weights, and $\widetilde R_k$ is the subvector of $R(x_k)$ formed by the components with indices $i\in {\cal{M}}_k$. With respect to the notation in Algorithm \ref{algo}, it holds $d=|{\cal{M}}_k|$.

\subsection{ {  \mq and \sistema  problems: inexact Gauss-Newton method with Jacobian sampling}}\label{secmodsistema}
In this section we consider {   problems (\ref{ls}) and  (\ref{sistema}) and specialize Algorithm 
\ref{algo_general} to the case where the Jacobian matrix is approximated by sampling and has the same dimension as $J(x_k)$; we will denote such random estimate as $\widetilde J_k\in \mathbb{R}^{m\times n}$.

For problem (\ref{ls}),  given $x_k$, we will use the  random model 
\begin{equation}\label{modelLS2}
\textstyle\widetilde m_k(s)=\frac{1}{2m}\|\widetilde J_k s+ R(x_k)\|^2,
\end{equation}
and let $g_k = \widetilde J_k^T R(x_k)$ be the corresponding estimate of the gradient  $\nabla f(x_k)$.
We remark that the model above has  row dimension $m$ while the residual in the  model (\ref{modelLS}) has $d\le m$ rows.

For problem (\ref{sistema}), given $x_k$, the stochastic counterpart of (\ref{model_nse}) is given by 
\begin{equation}\label{modelNSE}
    \widetilde m_k(s)=\frac 1 2 \|\widetilde J_k s+F(x_k)\|^2,
\end{equation}
 and $g_k = \widetilde J_k^TF(x_k)$ is  the corresponding estimate of the gradient  $\nabla f(x_k)$.

The inexact Newton step $s_k$ can be computed applying a Krylov solver with stopping criterion
\be \label{inexactgm2}
   \|\widetilde J_k^T  r_k\| \leq \eta_k \| \widetilde J_k ^T R(x_k)\| \, \, \mbox{ with }  \,\, r_k=\widetilde J_k s_k+  R(x_k),
\ee
for some  $\eta_k\in[0,\bar\eta), 0<\bar\eta<1$, in case of \mqp, and 
\be \label{inexactgm}
   \|\widetilde J_k^T  r_k\| \leq \eta_k \| \widetilde J_k ^T F(x_k)\| \, \, \mbox{ with }  \,\, r_k=\widetilde J_k s_k+  F(x_k),
\ee
for some  $\eta_k\in[0,\bar\eta), 0<\bar\eta<1$, in case of \sistemap.}

Algorithm \ref{algoNSE} describes the procedure sketched above.

\algo{}{Step 1 of Algorithm \ref{algo_general} for \mq and \sistema}{
\noindent
Given  $ x_k\in {\mathbb R}^n$.
\begin{description}
\item{Step 1.1} Choose $\eta_k\in {  [0,\bar \eta)}$. 
\item{Step 1.2} {  If the problem is of the form (\ref{ls}),
\\
\hspace*{25pt} Form $\widetilde J_k^T \in\RR^{m\times n}$, $g_k=\frac{1}{m}\widetilde J_k^T R(x_k)$.
\\
\hspace*{25pt} 
  {Apply a Krylov method to $\min_{s}\widetilde m_k(s)$ 
 with  $\widetilde m_k(s)$ given in (\ref{modelLS2}) 
 \\
 \hspace*{25pt} and compute $s_k$ satisfying (\ref{inexactgm2})}.}
\\
\hspace*{15pt} Else\\
\hspace*{25pt} {   Form  $\widetilde J_k^T \in\RR^{n\times n}$, $g_k=\widetilde J_k^T F(x_k)$.
\\
\hspace*{25pt}   {Apply a Krylov method   to $\min_{s}\widetilde m_k(s)$ 
 with  $\widetilde m_k(s)$ given in  (\ref{modelNSE}) 
 \\
 \hspace*{25pt}   and compute $s_k$
satisfying (\ref{inexactgm})}.}
\end{description}
}\label{algoNSE}
\vskip 5pt
 {As for the Krylov methods,  it is recommended to use LSMR \cite{LSMR}, LSQR \cite{PS}; in case $\widetilde J_k$ is square and symmetric  MINRES-QLP \cite{minresqlp} is recommended\footnote{See https://web.stanford.edu/group/SOL/home\_software.html}. These methods
can compute a minimum-length solution of $\min_s \widetilde m_k(s)$}.
\begin{lem} \label{discesa}
Let $s_k, \widetilde J_k, g_k$ as in Algorithm \ref{algoNSE}.
Then
$
s_k^Tg_k\le 0.
$
\end{lem}
\begin{proof}
{  The claim follows easily proceeding as in Lemma  \ref{discesa_ls}.}
\end{proof}

To complete the description of Algorithm \ref{algoNSE}, we focus on the construction of $\widetilde J_k$. We apply sampling \minor{by} interpreting $J(x_k)$ as a sum of matrices and consider two different approximations; in one case $J(x_k)$ is the sum of sparse and rank-1 matrices and  we form a sparse approximation, in the other case $J(x_k)$ is the sum of Jacobians, as in finite sum-minimization, and we construct a standard  batch approximation \cite{bbn_2017,bbn}.

{  The use of a sparse approximation of the Jacobian that still retains significant information in the Jacobian, has potential advantages: reduction of both
the storage requirement and   the cost of matrix-vector computations needed in the Krylov iterative solver. It  can be effective when the Jacobians  are dense and  too large to handle, when small elements can be zeroed out though retaining most of the information, and when redundant information is present such as in data analysis, see e.g.,  \cite{GR, HCZ, tropp}. Sparsification can be attempted   replacing the Jacobian by a sparse matrix that is close in  terms of spectral-norm distance. }
This task can be performed randomly selecting a small number of entries from the original matrix.
Let \minor{us} denote $E_{ij}$ the matrix that has the element in position $(i,j)$ equal to 1 and zeros otherwise, and   denote $J(x_k)_{(i,j)}$   the $(i,j)$ entry of $J(x_k)$, then
$$J(x_k) = {  \sum_{i=1}^m}\sum_{j=1}^n J(x_k)_{(i,j)}E_{ij}.$$
 Following \cite[\S 6.3]{tropp} and letting ${\cal{M}}_k$ be a subset if indices $(i,j)$ of all entries, we can generate a random approximation $\widetilde J_k$  by sampling as 
 \begin{equation}\label{Jsparse}
\widetilde J_k =   \frac{1}{|{\cal{M}}_k|}\sum_{(i,j)\in {\cal{M}}_k} \frac{1}{p_{i,j}^k}J(x_k)_{(i,j)}E_{ij}.
 \end{equation}
Matrix $\widetilde J_k$ is an unbiased estimator of $J(x_k)$.

The probability distribution can be  assumed  uniform, ${   p_{ij}^k=\frac {1}{mn}}$  $\forall i, j$,  or of the form associated to the so-called importance sampling  \cite[\S 6.3.3]{tropp}.
Given an accuracy requirement  $\rho>0$ and $\delta_J\in(0,1)$, it holds
$$\mathbb P(\|J(x_k)-\widetilde J_k\|\leq\rho)\geq 1-\delta_J$$
whenever the size of the sample ${\cal{M}}_k$  is sufficiently large according to Theorem \ref{Bernth}.

\REV{As a second type of sampling, we 
suppose that $J(x_k) $ is the  average  of  $N$ matrices,  say  $J(x_k) =  \frac{1}{N}\sum_{i=1}^N \Theta_i(x_k)$,  and  let   $p_i^k=\frac{1}{N}$, $i=1, \ldots, N$, denote the uniform probability distribution associated to matrices $\Theta_i(x_k)$.
Given a set ${\cal{M}}_k$ generated by randomly and uniformly sampling the set of indices $\{1,\dots,N\}$, it holds 
\begin{equation}\label{Jsomme}
\widetilde J_k = \frac{1}{|{\cal{M}}_k|}\sum_{i\in {\cal{M}}_k} \Theta_i(x_k).
\end{equation}
This matrix is an unbiased estimator of $J(x_k)$ and 
the sample size $|{\cal{M}}_k |$ which provides (\ref{Jsomme}) is again provided by  Theorem \ref{Bernth}.
}

\section{Iteration complexity for first-order optimality}\label{glob}
The algorithms introduced in the previous sections generate a stochastic process. 
Following \cite{cartis}, we denote $ \mathcal T_k$ the random
step size parameter, $\cS_k$   the random search direction, $X_k$   the random  iterate, and $ \mathcal J_k$ the random matrix used either  in (\ref{modelLS}) or in  (\ref{modelNSE}). 
Given $\omega_k$ from a proper probability space, we denote the realizations of the random variables above as
$t_k = \mathcal T_k(\omega_k)$, $s_k = \cS_k(\omega_k)$, $x_k=X_k(\omega_k)$, and $\widetilde J_k=\mathcal{J}_k(\omega_k)$. For brevity we will omit  $\omega_k$ in the following.
Given $x_{k}$ and $t_{k}$, the Jacobian estimator $\mathcal J_k$ generates the gradient estimator $\mathcal{G}_k$
of $f$. We use $\mathcal{F}_{k-1} = \sigma(\mathcal J_0,\ldots, \mathcal J_{k-1})$
to denote the $\sigma$-algebra generated by $\mathcal J_0, \ldots, \mathcal J_{k-1}$, 
up to the beginning  of iteration $k$.

In this section we study the  properties of the presented  algorithms and  provide the expected  number of iterations required to reach an $\epsilon$-approximate first-order optimality point, i.e., a point $x_k$ such that $\|\nabla f(x_k)\|\le \epsilon$ for some positive scalar $\epsilon$. 

Our analysis first derives technical results on the relationship between the trial step $s_k$ and the stochastic gradient $g_k$, then analyzes the occurrence of successful iterations, and finally obtains the expected iteration complexity bound relying on the  framework provided in \cite{cartis}. We start by making the following basic assumption. 
\begin{ass}\label{soluz}(Existence of a solution)
There exists a solution of problem (\ref{ls}). Problem  (\ref{problem1})  admits a zero residual solution.
\end{ass}
Moreover, for any realization of the algorithm, given the Jacobian $J(x_k)$ of the residual functions at $x_k$, we denote its singular value decomposition as 
$J(x_k) = U_k\Sigma_k V_k^T$, where $U_k, V_k$ are orthonormal, $\Sigma_k=\diag(\sigma_{k,1},\ldots,\sigma_{k,n})$,
$\sigma_{k,1}\geq\ldots\geq\sigma_{k,r}> \sigma_{k, r+1}=\ldots=\sigma_{k,n}=0$, with $r$ being the rank of the matrix; concerning matrix dimensions, it holds $U_k\in \mathbb{R}^{m\times m}$, $V_k\in \mathbb{R}^{n\times n}$ for problem (\ref{ls}),  
$U_k, V_k\in \mathbb{R}^{n\times n}$ for (\ref{problem1}).
The rank retaining factorization is denoted as
 \begin{equation}\label{svdrr}
J(x_k) = U_{k,r}\Sigma_{k,r} V_{k,r}^T,
\end{equation}
where $U_{k,r}, V_{k,r}$ denote the first $r$ columns of $U_k, V_k$ and
 $\Sigma_{k,r} =\diag(\sigma_{k,1},\ldots,\sigma_{k,r})$. For matrix $\widetilde J_k$ we denote its rank with $\widetilde r$, its singular values with $\widetilde \sigma_{k,i}$  and let 
$ \widetilde J_k = \widetilde U_k\widetilde \Sigma_k \widetilde V_k^T$ be the singular value decomposition 
and 
\begin{equation}\label{svdrrtilde}
\widetilde J_k = \widetilde U_{k,\widetilde r}\widetilde \Sigma_{k,\widetilde r} \widetilde V_{k,\widetilde r}^T
\end{equation}
be the rank retaining factorization.
\subsection{Analysis of the trial step} 
We establish bounds on the trial step $s_k$ that are necessary to characterize successful iterations and consequently the generated sequence $\{x_k\}$.
These bounds hold whenever 
the nonzero eigenvalues of $\widetilde J_k^T\widetilde J_k$
are uniformly bounded from below and above for some sufficiently small $\sigma_{\min}$ and for some sufficiently large $\sigma_{\max}$.
Then, let us introduce the following event.

\begin{defi}{\rm (Spectral properties of $\widetilde J_k$)}\label{def:Jacc}
Let $ \mathcal J_k$  be generated in either  Algorithm \ref{algo}
or Algorithm \ref{algoNSE}, 
and ${\cal{E}}_k$   be the event
$$
{\cal{E}}_k =\1\left (\sigma_{\min}\le \sigma_i(  \mathcal J_k^T \mathcal J_k) \le  \sigma_{\max},\quad i=1,\ldots , \widetilde {\cal R} 
\right)
$$
where $\1$ denotes the indicator function of an event,  $\widetilde {\cal R}$ is  the random variable whose realization is $\widetilde r$ in (\ref{svdrrtilde}) and  $0< \sigma_{\min}\le \sigma_{\max}$. 
\end{defi}
\noindent 
In the following Lemma we provide conditions that ensure that  ${\cal{E}}_k=1$.
\begin{lem}\label{svd_tildej}
Let  $ x_k $ be given, $\widetilde J_k$ be generated by either Algorithm \ref{algo} or \ref{algoNSE} and  $r$ and $\widetilde r$ be the rank of $J(x_k)$ and $\widetilde J_k$, respectively. Assume that $ \Sigma_{k,  r}$ defined in (\ref{svdrr}) satisfies
$$
2\sigma_{\min} I_{ r}  \preceq   \Sigma_{k,  r}^2 \preceq \frac{\sigma_{\max}}{2} I_{ r},
$$
with $\sigma_{\min}, \sigma_{\max}$ as in Definition \ref{def:Jacc}. 
\begin{description}
\item{i)}
Consider the \mq problem and Algorithm \ref{algo}. Let $\minor{\overline{J}}_k\in \mathbb{R}^{m\times n}$ be the matrix \REV{given in \eqref{barJ-NLS}.} 
Then ${\cal{E}}_k=1$  {if $\|J(x_k)-\overline{J}_k\|\leq \min\{ \sigma_{k, \widetilde r}-\sqrt{\sigma_{\min}}, \sqrt{\sigma_{\max}}-\sigma_{k, 1}\}.$}

\item{ii)} Consider the \REV{\mq  and \sistema problems} and Algorithm \ref{algoNSE}.
Then ${\cal{E}}_k=1$ if $\widetilde r\le r$ and  {if $\|J(x_k)-\widetilde{J}_k\|\leq \min\{ \sigma_{k, \widetilde r}-\sqrt{\sigma_{\min}}, \sqrt{\sigma_{\max}}-\sigma_{k, 1}\}.$}
 \end{description}
\end{lem}
\begin{proof}
$i)$ 
Let $\overline{J}_k$ and $\widetilde J_k$ be the matrices introduced in \S \ref{secmodmq}. The interlacing property of singular values decomposition gives that the rank of $\minor{\overline{J}}_k$ is at most $r$, \cite[Theorem 7.3.9]{HJ}. 
Further, letting $\minor{\Delta_k}=J(x_k)-\minor{\overline{J}}_k$ and  $\sigma_{k,i},\,  \bar \sigma_{k,i}$, $i=1\ldots,n$, be
the singular vales of $J(x_k),\, \minor{\overline{J}}_k$ respectively, 
we know that $|\sigma_{k,i}-\bar \sigma_{k,i}|\le \|\minor{\Delta_k}\|$, $\forall i$ \cite[Corollary 7.3.8]{HJ}.  
Thus, from the assumption on $\Sigma_{k,r}$ it follows that   the singular values $\bar \sigma_{k,i}$, $i=1, \ldots, r$ are uniformly bounded from below and above by  $\sigma_{\min}$ and   $\sigma_{\max}$ respectively when $\|\minor{\Delta_k}\|\le \min\{ \sigma_{k, \widetilde r}-\sqrt{\sigma_{\min}}, \sqrt{\sigma_{\max}}-\sigma_{k, 1}\}$. Since the singular values of $\minor{\overline{J}}_k$ and $\widetilde J_k$ are  equal, the thesis follows.

$ii)$
Let  $\widetilde J_k$ be the matrix introduced in \S \ref{secmodsistema}. 
Letting $\minor{\Delta_k}=J(x_k)-\widetilde J_k$ and  $\sigma_{k,i},\,  \widetilde \sigma_{k,i}$, $i=1\ldots,n$, be
the singular values of $J(x_k),\, \widetilde J_k$ respectively,
we know that $|\sigma_{k,i}-\widetilde \sigma_{k,i}|\le \|\minor{\Delta_k}\|$, $\forall i$ \cite[Corollary 7.3.8]{HJ}. 
Thus, from the assumption on $\Sigma_{k,r}$ it follows that   ${\cal{E}}_k=1$ whenever $\|\minor{\Delta_k}\|\le \min\{ \sigma_{k, \widetilde r}-\sqrt{\sigma_{\min}}, \sqrt{\sigma_{\max}}-\sigma_{k, 1}\}$.
\end{proof}  
The following lemma establishes useful technical results on $s_k$.
\begin{lem} \label{lemmafund}
Let $x_k$ be given, {   $\widetilde J_k$, $s_k, g_k$}  be generated  either in Algorithm \ref{algo} or in Algorithm \ref{algoNSE}.  
Suppose that  ${\cal{E}}_k=1$ 
 and that 
\be\label{Vs}
\|\widetilde V_{k,\widetilde r}^Ts_k\|\ge \mu\|s_k\|,
\ee
with $\widetilde V_{k,\widetilde r}$ defined in (\ref{svdrrtilde}).
{  Then,  
\begin{equation}\label{kconstant} 
\kappa_2\|g_k\| \leq \| s_k \| \leq \kappa_1 \|g_k\|, \quad -g_k^Ts_k \geq \frac{1}{\kappa_1}\| s_k \|^2, 
\end{equation}
for some positive constants $\kappa_1, \kappa_2$ independent of $k$.
}
\end{lem}
\begin{proof}
First consider problem \eqref{ls} and  Algorithm \ref{algo}. 
By (\ref{stoplsqr}) we have 
$g_k={  \frac{1}{m}} \widetilde J_k^T \minor{\widetilde{R}_k} = {  \frac{1}{m}}\widetilde J_k^T(-\widetilde J_k s_k+r_k)$ and 
$$
 \|g_k\| \le {  \frac{1}{m}}( \sigma_{\max}\|s_k\|
 +  \eta_k \|\widetilde J_k^T \minor{\widetilde{R}_k}\|)\le {  \frac{\sigma_{\max}}{m}}\|s_k\|+ \bar \eta \|g_k\|.
$$
Thus,  the leftmost in {   $\frac{(1-\bar \eta)m}{\sigma_{\max}}$}.
Further, 
(\ref{inmq}) and (\ref{Vs}) imply
\begin{eqnarray*}
 -  s_k^T g_k  & =  & {  \frac{1}{m}}  s_k^T \widetilde V_{k,r} \widetilde  \Sigma_{k,r}^2 \widetilde V_{k,r}^T s_k
\geq {  \frac{\sigma_{\min}}{m}} \|\widetilde V_{k,r}^Ts_k\|^2   \geq {  \frac{\sigma_{\min}}{m}} \mu^2\|s_k\|^2,
\end{eqnarray*}
i.e.,  the second inequality and the rightmost part of (\ref{kconstant}) hold with $\kappa_1={  \frac{m}{\sigma_{\min}\mu^2}}$.

{  The same results holds for problem \eqref{ls}  and Algorithm \ref{algoNSE} using (\ref{inexactgm2}) and 
$g_k=\frac{1}{m} \widetilde J_k^T {R}_k = \frac{1}{m}\widetilde J_k^T(-\widetilde J_k s_k+r_k)$}.

{  The claim for Algorithm \ref{algoNSE}  when applied to problem \eqref{sistema} follows using the same arguments as above, and the constants $\kappa_1$ and $\kappa_2$ take the form $\kappa_1=\frac{1}{\sigma_{\min}\mu^2}$ and 
$\kappa_2=\frac{1-\bar \eta}{\sigma_{\max}}$, respectively.}

\end{proof}
{  
\subsection{Regularized models}
The conditions  (\ref{Vs}) assumed in our analysis can be enforced using Levenberg-Marquardt models, i.e., regularizing the random models $\widetilde m_k(s)$ with a Tikhonov regularization and parameter 
$\lambda_k\in [\lambda_{\min},  \lambda_{\max}]$ for some $\lambda_{\min}$ and $\lambda_{\max}$ positive and small. E.g.,  in case of problem (\ref{problem1}) if
$$
\widetilde m_k(s)+\frac 1 2 \lambda_k\|s\|^2= \frac 1 2 \left\| \left (
\begin{array}{c}
\widetilde J_k\\ \sqrt{\lambda_k} I
\end{array}
\right ) s+ \left (\begin{array}{c}
F(x_k)\\ 0
\end{array}
\right )
\right \|^2,
$$
then 
$$
-s_k^Tg_k=s_k^T (\widetilde J_k^T  \widetilde J_k +\lambda_k I)s_k\ge 
\lambda_{\min}\|s_k\|^2.
$$
Hence, \eqref{kconstant} holds with $\kappa_1 = \frac{1}{\lambda_{\min}}$ and $\kappa_2 = \frac{1-\bar \eta}{\sigma_{\max}+\lambda_{\max}}$ and assumption  (\ref{Vs}) is no longer necessary.
Moreover, Definition \ref{def:Jacc} becomes
$$
{\cal{E}}_k ={\mathbbm 1}\left (\max_{i} \{\sigma_i(  \mathcal J_k^T \mathcal J_k+\lambda_k I)\} \le  \sigma_{\max}\right),
$$
since $\sigma_{\min}\ge \lambda_{\min}$, and proceeding as in Lemma \ref{svd_tildej} we get that $\mathcal{E}_k=1$ if $\|J(x_k)-\widetilde{J}_k\|\leq \sqrt{\sigma_{\max}}-\sigma_{k, 1}.$ 
}
\subsection{Fulfillment of the Armijo  condition} \label{sec3.2}
The study of the stochastic sequence $\{x_k\}$ depends on characterizing successful iterations and requires to assume accurate derivatives with fixed probability. 
In case of {  row-compression of the Jacobian}\minor{,} we assume that the stochastic gradient is sufficiently accurate with respect to  $\nabla f(x_k)$ in probability.
\begin{ass}
{\rm (gradient estimate, {   row compression of the Jacobian})}\label{gaccLM} 
Let $\alpha$ be a positive constant and {  consider Algorithm \ref{algo}}. The estimator $\mathcal J_k$ is $(1-\delta_g)$-probabilistically sufficient\minor{ly} accurate
in the sense that the indicator variable 
\begin{equation}\label{accgrad}
\textstyle{\cal{I}}_k = \1 \left ( \|  \nabla f(X_k)-{\cal{G}}_k \| \leq \alpha t_k \|{\cal{G}}_k\|  \right)
\end{equation}
satisfies the  submartingale condition
\be\label{probLM}
 \mathbb P\left( {\cal{I}}_k=1 | \mathcal{F}_{k-1}  \right) \geq 1 - \delta_g, \quad  \delta_g\in(0,1).
\ee
\end{ass}
This requirement can be satisfied approximating $\nabla f(x_k)$
by sampling as described in \S \ref{secmodmq}.  In this regard, note that the cardinality $|{\cal{M}}_k| $ depends on $\rho=\alpha t_k\|g_k\|$ in (\ref{accg}) with $\|g_k\|$ given in  (\ref{gk_mq}) but $g_k$ is unknown. In practice, one can enforce condition \eqref{probLM} proceeding as in \cite[Algorithm 4.1]{bg}.

In the {  remaining cases}, the Jacobian is supposed to be probabilistically accurate. 
\begin{ass}{\rm ({  Jacobian estimate})}\label{gacc} 
Let $\alpha$ be a positive constant and {  consider Algorithm \ref{algoNSE}}. The estimator $\mathcal J_k$ is $(1-\delta_J)$-probabilistically sufficiently accurate
in the sense that the indicator variable 
\begin{equation}\label{Jacc}
{\cal{I}}_k = \1 \left( \| J(X_k)- \mathcal J_k  \| \leq  \alpha t_k \right)
\end{equation}
satisfies the  submartingale condition
\be\label{deltaJ}
 \mathbb P\left( {\cal{I}}_k=1 | \mathcal{F}_{k-1}  \right) \geq 1 - \delta_J, \quad  \delta_J\in(0,1).
\ee
\end{ass}
\noindent
This accuracy requirement above can be fulfilled proceeding as in \S \ref{secmodsistema}.

Now we introduce the case where ${\cal{E}}_k={\cal{I}}_k=1$ holds and denote such occurrence as a {\em true} iteration.
\begin{defi}{\rm (True iteration)}
Iteration $k$ is true when  ${\cal{E}}_k{\cal{I}}_k=1$.
\end{defi}

\vskip 5pt
For true iterations a relevant  bound on $\nabla f(x_k)$ holds.
\vskip 5pt
\begin{lem} Consider any realization $x_k$  of Algorithm \ref{algo_general} and suppose that iteration $k$ is true. Suppose that Assumptions  \ref{gaccLM}, \ref{gacc} hold.
\begin{description}
 \item{i)} Consider the Algorithm \ref{algo}. It   holds
\be \label{truecons}
\| \nabla f(x_k) \| \leq \left( 1+ \beta t_k\right) \| g_k \|,
\ee
for some positive scalar $\beta$.
\item{ii)} Consider the Algorithm \ref{algoNSE}. Suppose that  {  
\be \label{up_boundF} 
{   
\| \widetilde J_k^T R(x_k)\|\ge \frac{1}{\kappa_3}  \|R(x_k)\|, \quad 
\|  \widetilde J_k^T F(x_k)|\ge \frac{1}{\kappa_3}  \|F(x_k)\|,}
\ee 
for problems \mq and \sistema respectively, and  some positive scalar $\kappa_3$.}
Then (\ref{truecons}) holds for some positive scalar $\beta$.
\end{description}
\end{lem}
\begin{proof}
$i)$ 
The claim follows trivially by Assumption \ref{gaccLM} with $\beta=\alpha$.
\\
$ii)$ Using \eqref{Jacc}, for problem \mq we obtain
\begin{eqnarray*}
\|\nabla f(x_k)\|  &\le& \| \nabla f(x_k) -g_k \| +\|g_k\|\\
&=& {  \frac{1}{m} }\| (J(x_k) -\widetilde J_k)^T R(x_k) \| +\|g_k\|\\
&\leq& {  \frac{1}{m} }\alpha t_k\|R(x_k)\|+\| g_k \|.
\end{eqnarray*}
Then by  (\ref{up_boundF})  and (\ref{kconstant})  we get 
{  $\| \nabla f(x_k) \| \leq \left(  1+ \alpha \kappa_3 t_k\right) \| g_k \|$
and the claim holds with $\beta= \alpha\kappa_3 $.} 
{  The same arguments can be applied for problem \sistemap.}
\end{proof}

Now we prove that if the iteration is true
and $t_k$ is small enough, the line-search condition is satisfied; namely, the iteration is successful. Thereafter we make the following assumption.
\vskip 5pt

\begin{ass}{\rm (gradient of $f$  Lipschitz-continuous)} \label{gradilipass} 
The gradient $\nabla f$ of $f$ is Lipschitz-continuous with constant $L$
\be \label{gradlip}
\| \nabla f(x) -\nabla f(y)\| \leq L\|x-y \| \text{ for all } x,y\in\RR^n.
\ee
\end{ass}

\begin{lem}\label{lowertk} 
Consider any realization of Algorithm \ref{algo_general} and suppose that iteration $k$ is true.
Suppose that Assumptions \ref{gaccLM}, \ref{gacc} and  \ref{gradilipass} hold, that the assumptions of Lemma \ref{lemmafund} hold {  and that (\ref{up_boundF}) is satisfied}.
Then
there exists a positive scalar $\bar t$ independent of $k$ such that the iteration is successful  whenever $t_k \leq \bar{t}$.
\end{lem}

\begin{proof}
 Let $ k $ be an arbitrary iteration.  
Assumption \ref{gradilipass} implies, using the standard arguments for functions with bounded Hessians,
\begin{eqnarray*}
f(x_{k}+t_{k} s_{k})&=&
f(x_{k}) +  \int_{0}^{1} (\nabla f(x_{k}+y t_{k} s_{k}))^{T} ( t_{k}s_{k}) dy\\
& = & f(x_{k}) +  \int_{0}^{1} t_k(\nabla f(x_{k}+y t_{k} s_{k}) - \nabla f(x_k))^{T} s_{k} dy + t_k\nabla f(x_k)^Ts_k\\
& \leq & f(x_{k}) +  \int_{0}^{1} t_k\|\nabla f(x_{k}+y t_{k} s_{k}) - \nabla f(x_k)\| \| s_{k}\| dy + t_k\nabla f(x_k)^Ts_k\\
& \leq & f(x_{k}) +  \frac{L}{2}  t_{k}^2 \|s_{k}\|^2 + t_k\nabla f(x_k)^Ts_k.
\end{eqnarray*}
In the case of the {   Algorithm \ref{algo}},  inequality (\ref{kconstant}) and the definition of true iteration yield
\begin{eqnarray}
 f(x_k+t_ks_k) 
&\leq&  f(x_{k}) +  \frac{L}{2}  t_{k}^2 \|s_{k}\|^2 + t_k [\nabla f(x_k)- g_k]^Ts_k + t_k  g_k^Ts_k \nonumber\\
&\leq &  f(x_{k}) +  \frac{L}{2}  t_{k}^2 \|s_{k}\|^2 + \alpha t_k^2 \|g_k\|\|s_k\|+  t_kg_k^Ts_k \nonumber\\
&\leq&  f(x_{k}) +  \frac{L}{2}  t_{k}^2 \|s_{k}\|^2 + \frac{\alpha}{\kappa_2} t_k^2 \|s_k\|^2+  t_kg_k^Ts_k  \nonumber.
\end{eqnarray}
Using  (\ref{kconstant}) we have $-(1-c) g_k^Ts_k \geq(1-c)  \frac{1}{\kappa_1} \| s_k \|^2 $. Thus,
 if 
\begin{eqnarray*}
t_k \|s_k\|^2 \left( \frac{L}{2} +  \frac{\alpha}{\kappa_2}\right) \leq (1-c) \frac{1}{\kappa_1} \| s_k \|^2,
\end{eqnarray*}
then (\ref{armijo}) holds and the claim follows with 
$t_k \leq \bar{t}= \frac{2(1-c)\kappa_2}{(\kappa_2 L+2\alpha)\kappa_1}$.

Consider now {   Algorithm \ref{algoNSE} and  problem \eqref{problem1}}. By  the definition of true iteration we obtain
\begin{eqnarray}
 f(x_k+t_ks_k) &\leq&  f(x_{k}) +  \frac{L}{2}  t_{k}^2 \|s_{k}\|^2 + t_k \nabla f(x_k)^Ts_k \pm t_k  g_k^Ts_k  \label{disugdim}\\
&=&  f(x_{k}) +  \frac{L}{2}  t_{k}^2 \|s_{k}\|^2 + t_k \left (J_k^TF_k- \widetilde J_k^TF_k\right )^Ts_k + t_k  g_k^Ts_k \nonumber\\
&=&  f(x_{k}) +  \frac{L}{2}  t_{k}^2 \|s_{k}\|^2 + t_k F_k^T\left (J_k- \widetilde J_k\right ) s_k + t_k  g_k^Ts_k \nonumber\\
&\leq&  f(x_{k}) +  \frac{L}{2}  t_{k}^2 \|s_{k}\|^2 + \alpha t_k^2 \|F_k\|\|s_k\| + t_k  g_k^Ts_k \nonumber.
\end{eqnarray}
Then,  (\ref{up_boundF})  and  (\ref{kconstant}) yield
$$
 f(x_k+t_ks_k)
 \leq   f(x_{k}) +  \frac{L}{2}  t_{k}^2 \|s_{k}\|^2 + 
 {  \alpha t_k^2  \frac{\kappa_3}{\kappa_2}\|s_k\|^2} + t_k  g_k^Ts_k.
$$
Proceeding as above, the claim follows with  $t_k \leq \bar{t}=   \frac{2(1-c)\kappa_2}{(L\kappa_2 +2\alpha\kappa_3)\kappa_1}$.
{  The same reasoning applies to problem \eqref{ls}}.
\end{proof}

\subsection{Complexity analysis of the stochastic process}
{ In this section 
we provide a bound on the expected number of
iterations that our procedures take  in the worst case before they achieve a desired level of accuracy in the first-order optimality condition. The formal definition for such a number of iteration is given below.

\begin{defi}\label{ht}
Given some 
$\epsilon>0$, $N_\epsilon$ is the number of iterations required until $\| \nabla f(X_k) \| \leq \epsilon$ occurs for the first time. 
\end{defi}
The number of iterations $N_\epsilon$ is a random variable and it can be defined as the hitting time for our stochastic process. Indeed it has the property $\sigma(\1 \left(N_\epsilon > k\right )) \subset \mathcal{F}_{k-1}$.

Following  the notation introduced in Section \ref{sec3.2}
we let $X_k$, $k\ge 0$, be the random variable with realization $x_k=X_k(\omega_k)$ and  consider
the  following measure of progress towards optimality: 
\be \label{zk}
Z_k= f(X_0)-f(X_k).
\ee
Further, we let 
\be
Z_\epsilon = f(X_0) -f_{{\rm low}} = f(X_0), \label{ze}
\ee
be  an upper bound for $Z_k$ for any $k<N_\epsilon$, 
with $f_{{\rm low}}=0$ being the global lower bound of $f$.
 We denote with
$z_k=Z_k(\omega_k)$  a realization of the random quantity $Z_k$.
}

\begin{lem} \label{hrstrconv} 
Suppose that Assumptions  \ref{soluz}, \ref{gaccLM}, \ref{gacc} and \ref{gradilipass} hold. Suppose that  the assumptions of Lemma \ref{lemmafund} hold.
Suppose that iteration $k$ is true and consider 
any realization of Algorithm \ref{algo_general}. If the $k$-th iteration is true and successful, then 
\be
z_{k+1} \geq z_k + c \frac{\kappa_2^2}{\kappa_1} \frac{t_k}{(1+\beta t_{\max})^2} \| \nabla f(x_k) \|^2,
\ee
whenever $k<N_{\epsilon}$.
\begin{proof}
For every true and successful iteration, using (\ref{armijo}), (\ref{kconstant}) and (\ref{truecons}),   we have 
\begin{eqnarray*}
f(x_{k+1}) &\leq& f(x_k) + ct_k s_k^T g_k \\
&\leq & f(x_k) - ct_k \frac{\kappa_2^2}{\kappa_1}  \| g_k \|^2 \\
&\leq & f(x_k) - ct_k \frac{\kappa_2^2}{\kappa_1}  \frac{1}{(1+\beta t_k )^2} \| \nabla f(x_k) \|^2 \\
&\leq & f(x_k) - c \frac{\kappa_2^2}{\kappa_1} \frac{t_k}{(1+\beta t_{\max})^2} \| \nabla f(x_k) \|^2,
\end{eqnarray*}
and the last inequality holds since $t_k\leq t_{\max}$.
Now, changing the sign and adding $f(x_0)$ we conclude the proof.
\end{proof}
\end{lem}

\vskip 5pt

\begin{lem} \label{falsez}
Consider any realization of Algorithm \ref{algo}.
For every iteration that is false and successful, we have
$$
z_{k+1} > z_k.
$$
Moreover $z_{k+1}=z_k$ for any unsuccessful iteration.
\begin{proof}
For every false and successful iteration, using   (\ref{armijo}) and  $s_k^Tg_k\le 0$ (see Lemma \ref{discesa_ls}, Lemma \ref{discesa}) we have 
$$
{f(x_{k+1})  \leq  f(x_k) + ct_k s_k^T g_k \le   f(x_k). }
$$
Now, changing the sign and adding $f(x_0)$, the first part of the proof is completed. Finally for any unsuccessful iteration,  Step 2 of Algorithm \ref{algo_general}  gives $x_{k+1}=x_k$;  hence it holds $f(x_{k+1})=f(x_k)$   and  $z_{k+1}=z_k$.
\end{proof}
\end{lem}

To complete our analysis we need to assume that true iterations occur with some fixed probability.
\begin{ass}{\rm (probability of true iterations)}
\label{probtrue} 
There exists some $ \delta\in\left (0, \frac 1 2 \right)$ such that 
$$ \mathbb P\left( {\cal{I}}_k{\cal{E}}_k=1 | \mathcal{F}_{k-1}  \right) \geq 1 - \delta $$
\end{ass}
\noindent

Now we can state the main result on the expected value of the hitting time.

\begin{teo} \label{EN2}
Suppose that Assumptions  \ref{soluz}, \ref{gaccLM}, \ref{gacc} and \ref{gradilipass} and \ref{probtrue} hold. Suppose that  the assumptions  of Lemma \ref{lemmafund} hold.
Let $\bar t$ given in Lemma \ref{lowertk} and suppose $\bar t <t_0$.
Then   the stopping time $N_\epsilon$   of Algorithm\ref{algo_general} for the \mq and \sistema problems is bounded in expectation as follows
$$
\E[N_\epsilon] \leq \frac{2(1-\delta)}{(1-2\delta)^2}\left[  \frac{M}{\epsilon^2} + \log_{\tau}\frac{\bar t}{t_0} \right],
$$
with $M=\frac{(f(x_0)-f^*)(1+\beta t_{\max})^2\kappa_1}{c\kappa_2^2\bar t} $.
\begin{proof}
Let  
\be
h(t) = c \frac{\kappa_2^2}{\kappa_1}  \frac{t}{(1+\beta t_{\max})^2} \epsilon^2,
\ee
and note that $h(t)$ is non\minor{-}decreasing for $t\in [0,t_{\max}]$ and that $h(t)> 0$  for $t\in [0,t_{\max}]$.
For any realization $z_k$ of $Z_k$ in (\ref{zk}) of Algorithm \ref{algo} 
the following hold for all $k < N_\epsilon$:
\begin{enumerate}
\item[(i)] If iteration $k$ is true and successful, then $z_{k+1} \geq z_k +h(t_k)$ by Lemma \ref{hrstrconv}.
\item[(ii)] If $t_k \leq \bar t$ and iteration $k$ is true then iteration $k$ is also successful, which implies $t_{k+1} = \tau^{-1}t_k$ by Lemma  \ref{lowertk}.
\item[(iii)] $z_{k+1} \geq z_k $, for all $k$ ($z_{k+1} \geq z_k $ for all successful iterations by Lemma \ref{hrstrconv}
and \ref{falsez}); $z_{k+1} = z_k$ for any unsuccessful iteration $k$ by Lemma \ref{falsez}).
\end{enumerate}
Moreover,  our stochastic process $\{{\mathcal{T}_k}, Z_k\}$ obeys the expressions below.
By Lemma \ref{lowertk} and the definition of Algorithm \ref{algo} the update of the random variable $\mathcal T_k $ such that $t_k = \mathcal T_k(\omega_k)$  is
$$
\mathcal T_{k+1}=  \left\{
\begin{array}{ll}
\tau^{-1} \mathcal T_k & \mbox{ if } I_k=1 , \, \mathcal T_k \le \bar t \  \mbox{ (i.e., successful)}\\
\tau^{-1} \mathcal T_k & \mbox{ if the iteration is successful}, \,  I_k=0 , \, \mathcal T_k \le \bar t \\
\tau\, \mathcal T_k & \mbox{ if the iteration is unsuccessful}, \,  I_k=0 , \, \mathcal T_k \le \bar t \\
\tau^{-1} \mathcal T_k & \mbox{ if the iteration is successful},\mathcal T_k > \bar t \\
\tau\, \mathcal T_k & \mbox{ if the iteration is unsuccessful},   \mathcal T_k > \bar t \\
\end{array} 
\right.
$$
By  Lemma \ref{lowertk}  Lemma \ref{hrstrconv} and  Lemma \ref{falsez} the  random variable $Z_k $ obeys the expression
$$
Z_{k+1}\ge   \left\{
\begin{array}{ll}
Z_k +h(\mathcal{T}_k)& \mbox{ if } I_k=1 , \, \mathcal T_k \le \bar t \  \mbox{ (i.e., successful)}\\
Z_k & \mbox{ if the iteration is successful}, \,  I_k=0 , \, \mathcal T_k \le \bar t \\
Z_k & \mbox{ if the iteration is unsuccessful}, \,  I_k=0 , \, \mathcal T_k \le \bar t \\
Z_k +h(\mathcal{T}_k)  & \mbox{ if the iteration is successful}, I_k=1, \mathcal T_k > \bar t \\
Z_k& \mbox{ if the iteration is unsuccessful}, I_k=1, \mathcal T_k > \bar t \\
Z_k& \mbox{ if the iteration is unsuccessful}, I_k=0, \mathcal T_k > \bar t \\
\end{array} 
\right.
$$
Then Lemma 2.2--Lemma 2.7 and Theorem 2.1 in \cite{cartis} hold which gives the thesis along with the assumption $\delta<\frac 1 2$.
\end{proof}
\end{teo}

{  \subsection{Sample complexity}
In this  subsection, we provide  the total sample complexity of our procedures on the base of  the sample complexity analysis for adaptive optimization algorithms with stochastic oracles presented in  \cite{JSX2024}\footnote{
Our procedures fit in the framework given in
\cite[Algorithm 1]{JSX2024} with $f_k^0$  and $f_k^+$ in \cite{JSX2024} being the exact values  $f(x_k)$ and $f(x_k+t_k s_k)$,  $m_k$ in \cite{JSX2024} being $\tilde m_k$ given in either  (\ref{modelLS}) or  (\ref{modelLS2}), 
$H_k$  in \cite{JSX2024} being the random matrix $\widetilde J_k^T \widetilde J_k$, and the sufficient reduction  in \cite{JSX2024} represented by \eqref{armijo}.}.

Specifically, we suppose  ${\cal E}_k=1$,
and specialize the results in \cite{JSX2024} to our algorithms. We observe that the  total sample complexity of Algorithm \ref{algo}
represents the number of rows selected throughout Algorithm \ref{algo}, while the  total sample complexity of Algorithm \ref{algoNSE} represents either the number of entries in (\ref{Jsparse}) selected throughout the procedure or the number of terms in (\ref{Jsomme}) selected throughout the procedure.
Such sample complexities depend on  the prescribed accuracy requirements which involve, at iteration $k$, the steplength $t_k$ and the probability $\delta_g, \, \delta_J$. Following   \cite{JSX2024}, at each iteration we denote the sample complexity with parameter $t$ as $\sc(t)$
as for both the  algorithms $\delta_g$ and $\delta_J$  can be treated as a constant. Further, we let $\toc(k)=\sum_{j=1}^k \sc({\cal{T}}_j)$ be the random variable which denotes the total sample complexity of running our algorithms for $k$ iterations. 

Let use Theorem \ref{Bernth} for our random approximations(\ref{gk_mq}), (\ref{Jsparse}) and (\ref{Jsomme}),
and let $V^\dagger$ and $M$ be upper bounds, independent of $k$, for the per-sample second moments and the norm of the summands. 
\begin{itemize}
  
\item {\bf Algorithm \ref{algo}.} 
For carrying out the sample complexity analysis we consider the following accuracy requirement:
\begin{equation}\label{accgrad_new}
\textstyle{\cal{I}}_k = \1 \left ( \|  \nabla f(X_k)-{\cal{G}}_k \|\right) \leq \alpha t_k  \max \left\{{\varepsilon \psi}, \|{\cal{G}}_k\| \right\}\;\;\; \psi\in\left (0, \frac{1}{2\alpha t_{max}}\right)
\end{equation}
in place of \eqref{accgrad}. We stress that  
the  complexity result given in Theorem \ref{EN2} still holds under \eqref{accgrad_new} with $\bar t= \frac{2(1-c)\kappa_2}{\left(\kappa_2 L+2 \max\left \{ \alpha, \frac{1}{t_{\max}} \right \}\right)\kappa_1}$. Given  some positive $t$ at iteration $k$,  we have that 
\begin{eqnarray}
\sc(t) 
&\le&
\left( \frac{2V^\dagger}{( \alpha t \varepsilon\psi)^2 } + 
 \frac{4M}{3 \alpha t \varepsilon \psi }\right) \log\left(\frac{n+1}{\delta_g}\right). \label{comp1}
\end{eqnarray}
In fact, Theorem \ref{Bernth} implies that if the cardinality of the sample ${\cal M}_k$,   satisfies
 $$|{\cal M}_k|
\ge  \left( 
\frac{2V^\dagger}{(\alpha t_k \max \left\{\varepsilon \psi, \|g_k\|\right\})^2 } + 
\frac{4M}{3\alpha t_k \max \left\{ \varepsilon \psi, \|g_k\| \right\}} \right )\log\left(\frac{n+1}{\delta_g} \right) , 
 $$ 
 then the accuracy requirement given in \eqref{accgrad_new} holds with probability at least $1-\delta_g$.
 Then, a number of row of the order given in  the right-hand side of(\ref{comp1})  guarantees that Assumption 
 \ref{gaccLM} is satisfied.
 Taking into account that $\sc(t)$ is of the order of $\varepsilon^{-2}$,  
following the lines of \cite[Theorem 2, Theorem 4]{JSX2024} we can prove  that, if the scalar $\tau$ in Algorithm \ref{algo} is such that $\tau > (2\delta_g)^{\frac{1}{4}}$ then
$$
\E[\toc(\min\{N_\varepsilon, k\})] \leq 
O\left(\frac{V^\dagger}{\bar t^2\varepsilon^2} k\log_{\frac{1-\delta_g}{\delta_g}}(k)\,k^{2\log_{\frac{\delta_g}{1-\delta_g}}(\tau)}  \right)
= O\left(\frac{V^\dagger k}{\bar t^2 \varepsilon^2}\right).$$
Consequently, if Algorithm \ref{algoNSE} is run for $k=O\left(\frac{1}{\varepsilon^2}\right)$ iterations, the expected  total number of rows employed  is   
$O\left(\displaystyle \frac{V^\dagger}{\bar t^2\varepsilon^4}\right)$
where $\bar t$ is given in Lemma \ref{lowertk}.

\item {\bf Algorithm \ref{algoNSE}.} 
Using Theorem \ref{Bernth} we have that 
$$\sc(t) = \left( \frac{2V^\dagger}{(\alpha t)^2} + 
\frac{4M}{3\alpha t}\right) \log\left(\frac{n+1}{\delta_J}\right) .$$
Following again the lines of \cite[Theorem 2, Theorem 4]{JSX2024} we can prove that,  if the scalar $\tau$ in Algorithm   \ref{algoNSE} is such that $\tau > (2\delta_J)^{\frac{1}{4}}$ then   
{
$$
\E[\toc(\min\{N_\varepsilon, k\})] \leq  
O\left(\frac{V^\dagger}{\alpha^2\bar t^2}  k\log_{\frac{1-\delta_J}{\delta_J}}(k)\,k^{2\log_{\frac{\delta_J}{1-\delta_J}}(\tau)}  \right),
$$}
where $\bar t$ is given in Lemma \ref{lowertk}.
Therefore, if Algorithm \ref{algoNSE} is run for $k=O\left(\frac{1}{\varepsilon^2}\right)$ iterations  we can conclude that the expected total number of  entries  selected in (\ref{Jsparse})  or the total  number of terms employed in (\ref{Jsomme}) is  $O\left(\displaystyle \frac{V^\dagger}{\bar t^2 \varepsilon^2}\right).$

\end{itemize}
We conclude this section noting now that $k = \frac{1}{\varepsilon^2}$ is not  the number of iterations needed to reach the hitting time; in fact, Theorem \ref{EN2} provides a bound on the expected value of $N_{\varepsilon}$. On the other hand, by
$\E[N_\epsilon] \leq \frac{C}{\varepsilon^2}$ for some $C$ sufficiently large, it follows that the total sample complexity of Algorithms \ref{algo}  and \ref{algoNSE} is 
$O\left(\varepsilon^{-4}\right)$ and 
$O\left(\varepsilon^{-2}\right)$, respectively, in high probability \cite{JSX2024}.

}

\section{Numerical Results}
In this section, we study the numerical performance of the proposed methods 
and denote as {Algorithm \metls (Stochastic inexact Gauss-Newton method with Row Compression) the procedure for  the \mq problem, i.e., Algorithm \ref{algo_general} coupled with Algorithm \ref{algo}, 
and as Algorithm \REV{\metsist (Stochastic inexact Gauss-Newton method with Jacobian Sampling) the procedure for the  \mq and \sistema problems} consisting of Algorithm \ref{algo_general} coupled with Algorithm \ref{algoNSE}}.
We present the performance of \metls and \metsist algorithms and their full accuracy  counterparts, i.e., the algorithms  employing  exact Jacobians, indicated in the following as ``full''.
The parameters used in step 2 of Algorithm \ref{algo_general} are given by $c = 10^{-4}, t_{\max} = 1$ \REV{and $\tau=0.5$}. The values $\delta_g$ in (\ref{probLM}) and  $\delta_J$ in   (\ref{deltaJ})  
are equal to $0.4$.

\subsection{\REV{Nonlinear least-squares  with row-compression}}\label{exps_nls}
We consider the following least-squares problem
\begin{equation}\label{ls_loss}
\min_{x\in\R^n}f(x) = \|R(x)\|^2_2,\qquad R(x)_{(i)} = b_i-\frac{1}{1+e^{-x^T a_i}}\enspace i=1,\dots,m\end{equation}
where $a_i\in \R^n$, $b_i\in\{0,1\}$, $i=1,\dots,m$,  are the features vectors and the labels of the training set of a given binary classification problem.

We consider problem \eqref{ls_loss} for the \emph{gisette} dataset \cite{uci}, with $n=5000$ and $m=6000$. 
We also used a validation set of  1000 instances to evaluate the reliability of the classification model.
We run Algorithm \metls varying the parameter  $\alpha\in\{1,10,100\}$ 
in (\ref{accgrad})  and  constant forcing term $\eta_k = \eta = 10^{-1}$, $\forall k$. The matrix  $\widetilde{J}_k$ was generated by subsampling the rows of $J(x_k)$, with uniform probability, see \S 2.1. Concerning the cardinality of $\mathcal {M}_k$ we make use of  Theorem \ref{Bernth}, (\ref{accgrad}) and (\ref{probLM}), and choose $\mathcal {M}_k$ as follows:
\begin{equation}\label{sample_NLS}
    |\mathcal M_k|= \max \left\{0.01 m, \min\left \{ m_{\max},  2 \gamma\left(\frac{\|R(x_k)\|^2}{\rho_k^2}+\frac{2\|R(x_k)\|_{\infty}}{3\rho_k}\right)\log\left(\frac{ {n+1}}{\delta_g}\right)\right \} \right \},
    \end{equation}
with $\delta_g=0.4$,  $\rho_k= \alpha t_k \|g_{k-1}\|$, $m_{\max}\in [0.01m, m] $, $\gamma \in (0,1]$. Note that the accuracy request  \eqref{accgrad} is implicit and that 
$\rho_k$ in \eqref{sample_NLS}  employs $\|g_{k-1}\|$ instead of $\|g_k\|$ to make the evaluation of $ |\mathcal M_k|$ explicit with respect to the norm of the stochastic gradient.
We will report results varying $m_{\max}$ and $\gamma$, namely $m_{\max}\in \{ 0.75 m, \, m\}$ and $\gamma\in \{10^{-1},\,  1\}$. The choice $m_{\max}=m$ and $\gamma=1$ allows $\mathcal M_k$ to reach the full sample with the increase dictated by the Bernstein inequality, while the choice $m_{max}=m$ and $\gamma=10^{-1}$ retains the increase rate of the Bernstein inequality but  employing smaller  sample sizes. 
Clearly, the choice $m_{\max}=0.75m$ prevents the method from reaching the full sample. Note also that the size of  the sample is forced  to be at least  $1\%$ of 
$m$.

The initial guess $x_0 = 0\in\R^n$  was used and termination was declared  when either the number of full Jacobian evaluations is equal to 100, or the following stabilization condition holds for a number of iterations that corresponds to at least 5 full evaluations of the Jacobian
$$ |f(x_{k+1})-f(x_{k})|\leq \chi f(x_k)+\chi,$$
with $\chi = 10^{-3}$
\cite{sirtr_coap}. 

As for the computational cost, we assign cost $m$ to the evaluation of the residual vector $R(x)$, cost $n$ to the evaluation of one row the Jacobian, and cost  {     $2|\mathcal{M}_k|n$ } to the execution of one iteration of LSMR method \cite{LSMR} as \tblue{LSMR requires two matrix vector products (with matrices $\widetilde J_k$ and $\widetilde J_k^T $) at each iteration}, the resulting total cost was then scaled by the number of variables $n$. To summarize, the per-iteration cost of the method is given by
\REV{$$\frac{m}{n}+2\ell_k|\mathcal{M}_k|+|\mathcal{M}_k|,$$}
where $\ell_k$ is the number of inner iterations performed by the Krylov solver at $k$-th iteration.

\begin{figure}[h]
\centering
    \centering
        \includegraphics[width = 0.7\textwidth]{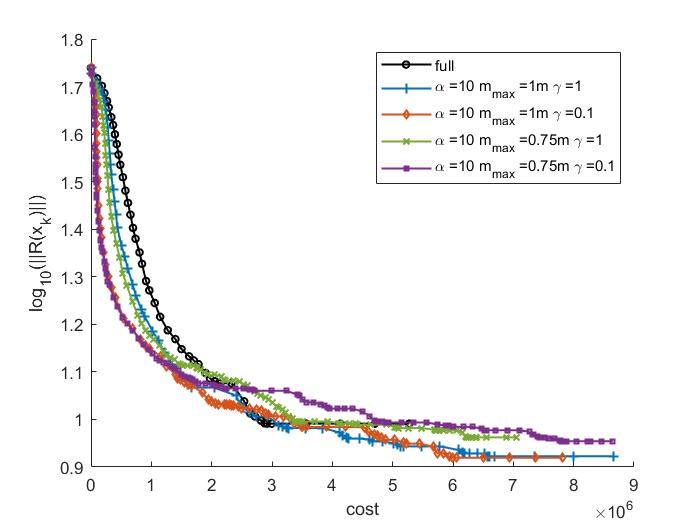}
    \caption{Algorithm \metlspunto, $\alpha=10$, varying $m_{\max}$ and $\gamma$. Median run in terms of cost: logarithmic norm of the residual versus computational cost.}\label{fig:ls_gisette}
\end{figure}

In Figure \ref{fig:ls_gisette} we report the objective function value, in logarithmic scale, versus the computational cost of 
the Algorithm \metls with $\alpha=10$ and varying  the values of $m_{\max}$ and $\gamma$. To account for the randomness of the Jacobian approximation, we run Algorithm \metls for each choice of the parameters 21 times and plot the results that correspond to the median run with respect to  the total computational cost at termination. We also plot the objective function value, in logarithmic scale, versus the computational cost of the full counterpart employing exact Jacobians.
With respect to this latter method, we see that the Algorithm \metls compares well in the initial stage of the convergence history and that attains smaller final values of $\|R\|$ when $m_{max}=m$.
Runs with $m_{max}=0.75m$, i.e., runs where full sample is not achieved, present a total computational cost at termination that is comparable to that of the algorithms with $m_{max} = m$ but larger values of the objective function at termination. 
 We also remark that the initial sample size is given by $600$ and $60$ for $\gamma =1$ and $\gamma = 0.1$ respectively, and that all the runs reach sample size $m_{max}$ when approaching termination.

\begin{figure}[h]
\centering
    \centering
        \includegraphics[width = 0.7\textwidth]{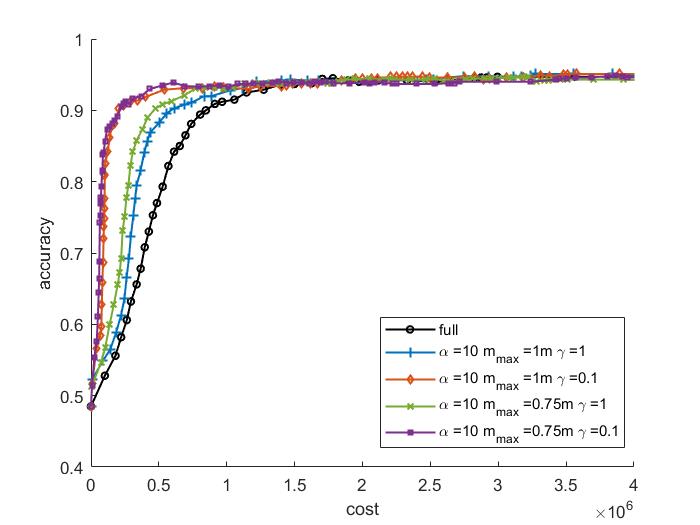}
    \caption{Algorithm  \metlspunto, $\alpha=10$, varying $m_{\max}$ and $\gamma$. Median run in terms of cost: accuracy versus computational cost.}\label{fig:ls_gisette_accuracy}
\end{figure}

In Figure \ref{fig:ls_gisette_accuracy} we report the accuracy, i.e., the percentage of entries of the validation set correctly classified versus the computational cost. In all runs, $0.94\%$ of the entries of the validation set is correctly classified and the figure shows that 
using row compression provides computational savings with respect to using the full Jacobian.  The figure displays the median run.
A similar behaviour is observed with $\alpha=1$. 

\subsection{\REV{Nonlinear least-squares and systems with Jacobian sampling}}
\REV{
In this section we discuss Algorithm \metsist applied to nonlinear least-squares problems and nonlinear systems. First, we consider random sparsification and adjust the sample size using either the Bernstein inequality and the importance sampling or fixed sample size with uniform sampling. Second, we address the solution  of one nonlinear systems representing  the first-order optimality conditions of a finite sum  {\em invex} objective function and use sampling for approximating the Jacobian.
}

\subsection{Importance sampling}
We apply Algorithm \metsist to \tblue{the} nonlinear system (\ref{sistema}) arising from the discretization of \tblue{an} integral equations. 
The system, named {\tt IE}, has equations of the form  
$$
F(x)_{(i)} = x_i+\frac{\left(1-h_i\right)}{2}\sum_{j=1}^i h_j\left(x_j+h_j+1\right)^3 + \frac{h_i}{2}\sum_{j=i+1}^n \left(1-h_j\right)\left(x_j+h_j+1\right)^2, $$
where $i=1, \ldots, n$, $n$ is the dimension of the system and  $h_j = j/(n+1)$, \cite{integr2}.

We set $n=5000$ and applied Algorithm \metsist using LSMR as the linear solver. The Jacobian approximation $\widetilde J_k$ was formed interpreting $J(x_k)$ as the sum of its diagonal part and its off-diagonal part and approximating the off-diagonal part of $J(x_k)$ by using (\ref{Jsparse}) with importance sampling, i.e., 
\begin{equation}\label{imp_s}
p_{ij}^k = \frac{1}{2}\left(\frac{|J(x_k)_{(i,j)}|^2}{\|J(x_k)\|_F^2} + \frac{|J(x_k)_{(i,j)}|}{\|J(x_k)\|_{\ell_1}} \right), \ \ i,j=1, n,
\end{equation}
with $\|J(x_k)\|_{\ell_1}= \sum_{i=1}^n\sum_{j=1}^n |J(x_k)_{(i,j)}|$,
and a number $|{\cal{M}}_k|$ of nonzero entries such that
\begin{equation}\label{pnounif}
n(n-1)\ge {|\cal{M}}_k|\geq \left(\frac{8\|J(x_k)\|_{\ell_1}}{3\alpha t_k}+\frac{4n \|J(x_k)\|_F^2 }{\alpha^2t_k^2}\right)\log\left(\frac{2n}{\delta_J}\right),
\end{equation}
which follows from Theorem \ref{Bernth} and the requirements (\ref{Jacc}), (\ref{deltaJ}).

The initial guess $x_0$ was drawn from the normal distribution $\mathcal{N}(0,1)$. 
Termination of Algorithm \metsist was declared  when $\|F(x_k)\|\leq 10^{-6}$.

\begin{figure}[h]
\centering
        \includegraphics[width = 0.8\textwidth]{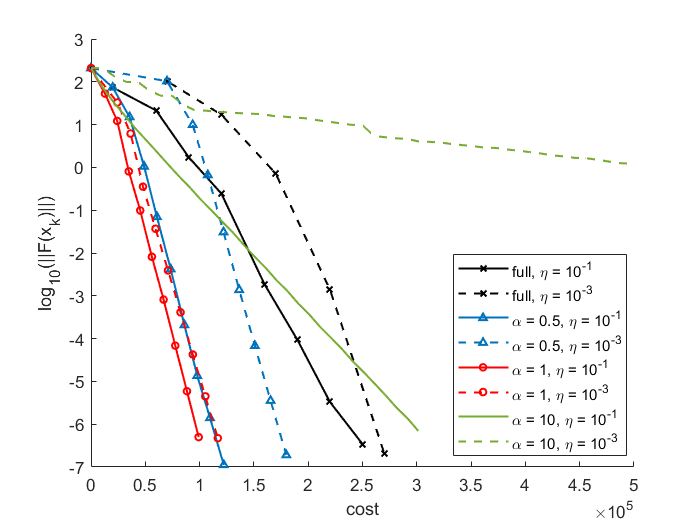}   
    \caption{Algorithm \metsist and Integral Equation {\tt IE}. Importance sampling.  Median run in terms of cost: logarithmic norm of the residual versus computational cost.}\label{fig:integr_more}
\end{figure}

Figure \ref{fig:integr_more} displays the results obtained in the solution of problem {\tt IE} testing three choices of the scalar $\alpha$ in (\ref{Jacc}),  \tblue{i.e., $\alpha\in\{0.5,1,10\}$}, and two choices of constant forcing terms, $\eta_k=\eta$, $\forall k$, $\eta\in\{10^{-3}, 10^{-1}\}$. 
We plot the computational cost and the norm of the residual $\|F(x_k)\|$ in logarithmic scale, on the $x$-axis and the $y$-axis respectively. The computational cost per iteration is evaluated as follows. 
Scaling again by the number of variables, we assign cost 1 to the evaluation of the vector $F\in\R^n$, cost $n$ to the evaluation of $J\in\R^{n\times n}$ as well as to the computation of the probabilities $\{p_{ij}^k\}_{i,j=1}^n$,  {and  $2(|{\cal{M}}_k|+n)/n$ to each iteration of LSMR.}  {   To summarize, the per-iteration cost of the method is given by 
$$
1+2n+2\ell_k\frac{|{\cal{M}}_k|+n}{n},
$$
where $\ell_k$ is the number of inner iterations performed by the Krylov solver at $k$-th iteration.}
To account for the randomness in the sparsification, each algorithm and parameter setting is run 11 times. In the plot we report the median run in terms of total computational cost at termination.

\tblue{ We first note the for $\alpha=0.5$ and $\alpha=1$, our algorithm is more convenient than Gauss-Newton method with exact Jacobian and that the best results are obtained using $\eta = 10^{-1}$. Setting $\alpha = 0.5$,  $\alpha=1$  and both  $\eta=10^{-1}$ and $\eta=10^{-3}$, all the runs with sparsification achieve the requested accuracy at a cost that is significantly smaller than that resulting from the use of the exact Jacobian. As expected, enlarging $\alpha$ reduces the accuracy of $\widetilde J_k$ and the value $\alpha=10$ deteriorates the performance of our algorithm. The most effective run corresponds to the use of $\alpha=1$ and $\eta = 10^{-1}$ and employed sparsified Jacobian with density varying between $0.03$ and  $0.28$. }

\begin{table}
\begin{center}

\begin{tabular}{ c|c|c|c|c } 
 \textbf{accuracy} & \textbf{cost} & \textbf{it} & \textbf{min cost} & \textbf{max cost}\\ 
  \hline
 full & 2.5001e+05 & 8 &      &    \\
 $\alpha = 0.5$ & 1.2226e+05	&9 	&1.0959e+05	&1.2227e+05	\\
 $\alpha = 1.0$ & 9.9123e+04	&10	&9.8863e+04	&9.9380e+04	\\
 
 $\alpha = 10$ & 3.0139e+05	&31	&3.0136e+05	&3.0142e+05 \\
\end{tabular}
\end{center}
\caption{\tblue{Algorithm \metsist and Integral equation {\tt IE}. Importance sampling, $\eta=10^{-1}$. Statistics on 
 multiple runs.}}\label{tab_es1_imp}
\end{table}
\tblue{Table \ref{tab_es1_imp} summarizes the results of multiple runs. For different values of $\alpha$, including the use of the exact Jacobian, it  displays the median cost (cost) of the runs, the number (it) of Gauss-Newton iterations performed in the run with median cost, the minimum and maximum cost of the multiple runs
(min cost, max cost). As expected, the number of iterations increases with $\alpha$ since the accuracy in the Jacobian approximation decreases. We also observe that the minimum and maximum value of the computational cost are close to the median values. }

\subsection{Uniform sampling }
\begin{figure}[h]
\centering
        \includegraphics[width = 0.8\textwidth]{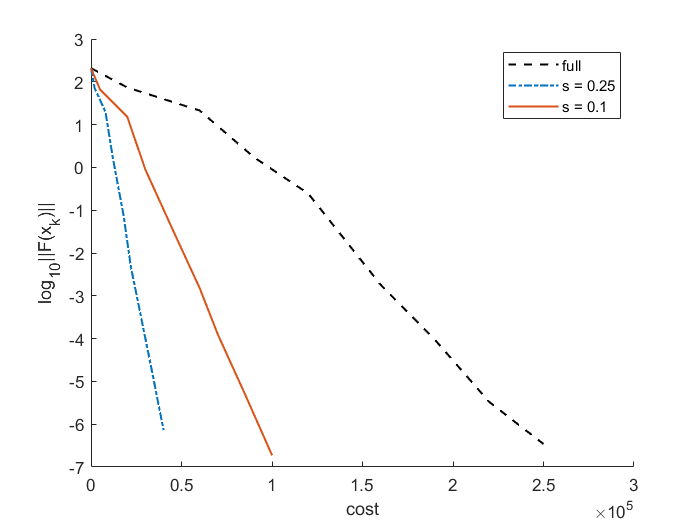}   
    \caption{Algorithm \metsist and Integral Equation {\tt IE}. Uniform sampling, density of the sparsified Jacobian equal to $s$. Median run in terms of cost: logarithmic norm of the residual versus computational cost.}\label{fig:integr_unif}
\end{figure}

\REV{The use of importance sampling  and probabilities $p_{ij}^k$ in  \eqref{imp_s} does not allow a matrix-free implementation and has an additional computational cost. We then investigated the behaviour of Algorithm  \metsist  employing Jacobian approximation obtained with uniform sampling and prefixed sample size without  relying on the Bernstein inequality.
}
\tblue{We solved {\tt IE} with $\widetilde J_k$ computed as follows.  Given a scalar $s\in(0,1)$, we let $\widetilde{J}_k$ be the matrix with the diagonal equal to the diagonal of $J(x_k)$ and $|{\cal{M}}_k| = sn^2-n$ off-diagonal elements uniformly sampled from $J(x_k)$. Consequently, $s$ represents the density of $\widetilde J_k$ and the computational per-iteration cost is given by
$ 2+sn+2\ell_k sn$ where $\ell_k$ is the number of inner iterations performed by the Krylov solver at $k$-th iteration.  In Figure \ref{fig:integr_unif} we plot the results  of the median run obtained setting $s=0.1$ and $s=0.25$,  and the  median run with  exact Jacobian. The forcing term used is constant,  $\eta_k =  10^{-1}$.
We can see from Figure \ref{fig:integr_unif} that, using sparsification, the required accuracy is achieved  with a smaller amount of computation than using the full Jacobian and that the best result is achieved with $s=0.25$, which takes 9 iterations to arrive at termination.}

\REV{A similar behaviour can be observed in the solution of nonlinear least-squares problems. We solved problem \eqref{ls_loss} using  at each iteration of \metsistp,  a matrix $\widetilde{J}_k$ with dimension  $m\times n$ and $|\mathcal{M}_k| = smn$ nonzero entries uniformly sampled from $J(x_k)$, with $s=0.1,\, 0.25, \, 0.5,\,  0.75$. }

\begin{figure}[h]
\centering
    \begin{subfigure}[b]{0.48\textwidth}
    \centering
        \includegraphics[width = \textwidth]{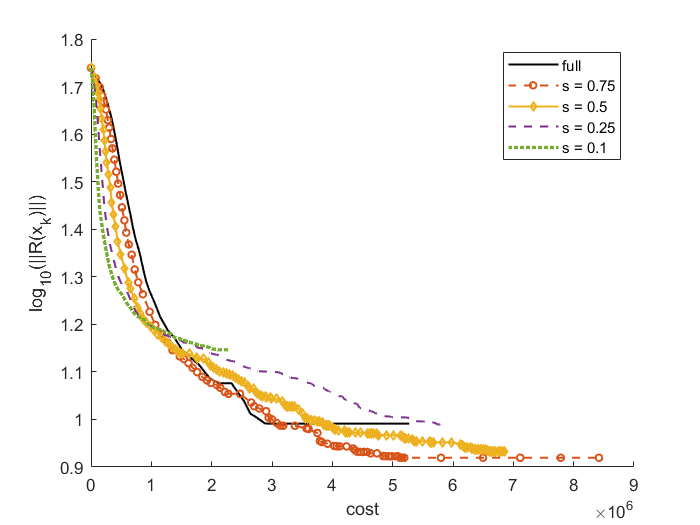}
 \caption{}
    \end{subfigure}
    \hfill
    \begin{subfigure}[b]{0.48\textwidth}
    \centering
        \includegraphics[width = \textwidth]{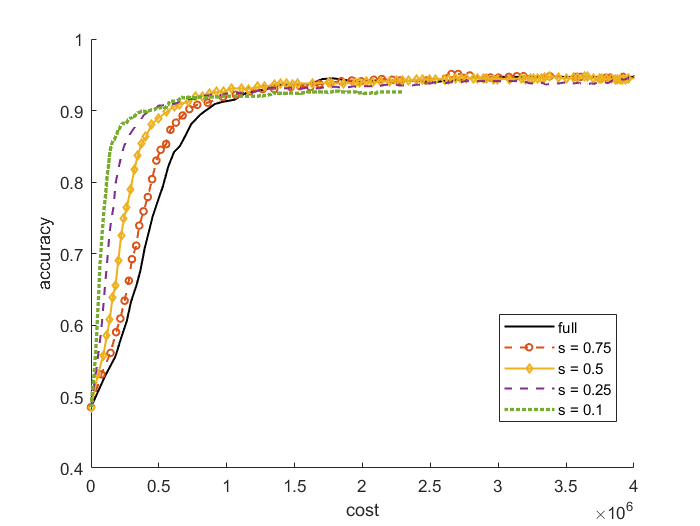}
     \caption{}
    \end{subfigure}
    \caption{Algorithm  \metsist with random sparsification. Median run in terms of cost. Computational cost versus logarithmic norm of the residual (a) and accuracy (b).   }\label{unif_ls}
    \end{figure}

\REV{
We run Algorithm \metsist  11 times for each choice of $s$, and in Figure \ref{unif_ls} we plot the results that correspond to the median run with respect to  the total computational cost at termination. 
As for the computational cost, we assign cost $m$ to the evaluation of the residual vector $R(x)$, cost $|{\cal{M}}_k|$ to the evaluation of the nonzero entries of the Jacobian, and cost    $2|\mathcal{M}_k|$  to the execution of one iteration of LSMR method, the resulting total cost was then scaled by the number of variables $n$.  To summarize, the per-iteration cost of the method is 
$$\frac{m+2\ell_k|\mathcal{M}_k|+|\mathcal{M}_k|}{n},$$
where $\ell_k$ is the number of inner iterations performed by the Krylov solver at $k$-th iteration.} 

\REV{
In Figure \ref{unif_ls} we plot the objective function value, in logarithmic scale,  and the accuracy, i.e., the percentage  of entries of the validation set correctly classified, versus the computational cost. 
We observe that the lowest values of the residual function are attained using $s=0.5$ and $s=0.75$. In terms of accuracy, the sparsified and full version of the algorithm reach the highest accuracy with comparable cost but the accuracy increases more rapidly in the case of sparsified Jacobians, especially for the values $s=0.1$, $s=0.25$ and $s=0.5$. }

We conclude this section with some comments on the potential savings resulting from random sparsification. The use of importance sampling  and probabilities $p_{ij}^k$ in  \eqref{imp_s} does not allow a matrix-free implementation and has  a cost that was taken into account in our measure for the computational burden. On the other hand,
sparsification by sampling provides saving in the Krylov solver and our experiments show that random models are overall advantageous. Finally, in the case of uniform probabilities,
forming $\widetilde J_k$ calls only for the evaluation of selected entries and \REV{ matrix-free implementations are possible. Our numerical experience shows that sparsification via  uniform probability distribution
can be effective  coupled with fixed sample sizes.}

\subsubsection{\REV{Softmax loss function: uniform subsampling}}
 Solving the  unconstrained optimization problem $\min_{x\in \R^n}\varphi(x),$
 with $\varphi:\R^n\rightarrow\R$ being an {\em invex} function is equivalent to solving the linear system of equations $F(x)= \nabla \varphi(x)=0$, see e.g., \cite{Newt-MR1}.
A binary classification problem performed via machine learning and the softmax cross-entropy convex loss function falls in such class and we solve such a problem in this section. 
The function $\varphi$ takes the form
\begin{equation}\label{softmax}
    \varphi(x) = \sum_{i=1}^N \varphi_i(x), \qquad  \varphi_i(x) = \log\left(1+e^{a_i^T x}\right) - \1(b_i=1)a_i^T x,
\end{equation}
where $\{a_i,b_i\}$, $i=1, \ldots, N$, is the dataset,  $a_i\in\R^n$, $b_i\in\{1,2\}$ and it is twice-continuously differentiable. In this section we report the  results of  the binary classification dataset \emph{a9a} \cite{uci},  with $n = 14 $ and $N = 30162$.

In the following, we apply Algorithm \metsist to the system $F(x) = \nabla\varphi(x)=0$; since  $J(x)$ is symmetric \REV{and square we used the iterative linear solver MINRES-QLP \cite{minresqlp}. Moreover, the approximate Jacobian $\widetilde J_k$ is formed by subsampling the sum $J(x_k)=\sum_{i=1}^N \nabla^2  \phi_i(x_k)$ as in (\ref{Jsomme}) and  using uniform probability distribution.} Using Theorem \ref{Bernth}, (\ref{Jacc}) and (\ref{deltaJ}), the rule for  $\mathcal{M}_k$  is 
\begin{equation}\label{Mk_Jsomme}
|{\cal{M}}_k | \geq \min \left\{N, \frac{4 \zeta(x_k)}{\alpha t_k}  \left(\frac{ {2}\zeta(x_k)}{\alpha t_k} +\frac{1}{3} \right) \log\left( \frac{2n}{\delta_J} \right) \right\}
\end{equation}
with $\max_{i\in\{ 1,\ldots,N \}}\| \nabla^2 \phi_i(x_k) \| \leq \zeta(x_k)$. It is know that such rule is expensive to apply as well as pessimistic, in the sense that it provides excessively  large values for $|{\cal{M}}_k |$. Hence  we applied (\ref{Mk_Jsomme})  setting $\zeta(x_k)=1, \forall k$, and 
\begin{eqnarray*}
& & |\mathcal{M}_k| = \max\{M_{\min}, \min\{N, \widehat{M}_k\} \},
\\& & M_{\min} = \xi N, \, \, \xi\in[0,1], \qquad 
\widehat{M}_k = \frac{4}{\alpha t_k}\left(\frac{1}{\alpha t_k}+\frac13\right)\log\left(\frac{2n}{\delta_J}\right).
\end{eqnarray*}
The parameter $\xi$ affects the value of $|\mathcal{M}_k|$, letting $\xi=1$ gives $|\mathcal{M}_k| = N$ at every iteration, i.e.,  $\widetilde J_k = J(x_k)$, $\forall k$; reducing $\xi$ may promote a reduction of $|\mathcal{M}_k|$.

We measure the computational cost at each iteration 
as follows. Let the cost of evaluating {     $\nabla  \varphi_i$ for any $i\in\{1,\dots,N\}$ be equal to $1$, thus evaluating $\nabla \varphi$ costs $N.$  } Each iteration of MINRES-QLP requires the computation of one Jacobian-vector product of the form $\widetilde{J}_k v = \sum_{j\in\mathcal{M}_k} \nabla^2\varphi_i(x_k) v, $ i.e., it requires $|\mathcal{M}_k|$ Hessian-vector products. Assuming that these products are computed with finite differences and taking into account that $\nabla \varphi_i(x_k)$ has already been computed at the beginning of the iteration to form $F(x_k)$, one MINRES-QLP iteration 
costs $|\mathcal{M}_k|.$ Consequently,  the $k$-th iteration of our algorithm  costs $N+|\mathcal{M}_k|\ell_k$ with $\ell_k$ being the number of MINRES-QLP iterations.

The  {   experimental results presented in the previous sections} indicate that the  computational cost of our solver depends on: the number of  nonlinear iteration performed; the cardinality of $\mathcal{{M}}_k$; the forcing terms $\{\eta_k\}$.
We investigate the choice of $\xi$ and $\eta_k$'s  by  applying Algorithm \metsist with $\xi\in\{1,10^{-1},10^{-2}, 10^{-3}, 0\},$ and constant forcing term $\eta_k = \eta, \forall k$, \tblue{$\eta \in\{ 10^{-1},10^{-3}\}$}.  
Further, we set $\alpha = 1$ in (\ref{Jacc}). 
The initial guess is $x_0 = 0$ and the stopping criterion   for the
algorithm is $\|F(x_k)\|\leq 10^{-3}.$ 
To account for the randomness in the method, for  each setting of $\xi$ and  $\eta$ out of the fifteen considered, the algorithm is run 21 times.

\begin{figure}[h]
\centering
    \begin{subfigure}[b]{0.48\textwidth}
    \centering
        \includegraphics[width = \textwidth]{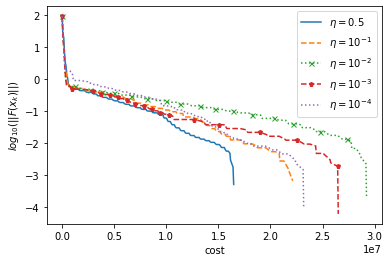}
 \caption{$\xi = 1$ (full sample)}
    \end{subfigure}
    \hfill
    \begin{subfigure}[b]{0.48\textwidth}
    \centering
        \includegraphics[width = \textwidth]{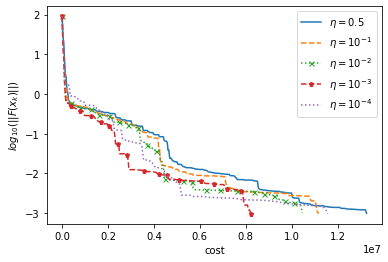}
     \caption{$\xi = 10^{-1}$}
    \end{subfigure}
    
    \begin{subfigure}[b]{0.48\textwidth}
    \centering
        \includegraphics[width = \textwidth]{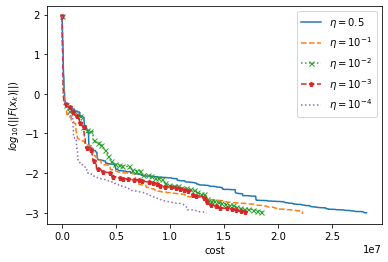}
 \caption{$\xi = 10^{-2}$}
    \end{subfigure}
    \begin{subfigure}[b]{0.48\textwidth}
    \centering
        \includegraphics[width = \textwidth]{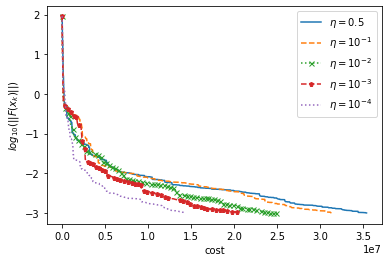}
 \caption{$\xi = 10^{-3}$}
    \end{subfigure}

    \begin{subfigure}[b]{0.48\textwidth}
    \centering
        \includegraphics[width = \textwidth]{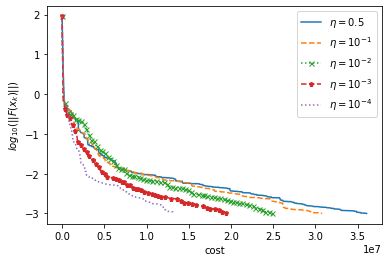}
 \caption{$\xi = 0$}
    \end{subfigure}
   
    \caption{\tblue{Algorithm \metsist and binary classification for the dataset  \emph{a9a}. Exact Jacobians ($\xi=1$)  and   sparsified Jacobians  with varying $\xi$ and $\eta$.   Median run in terms of cost: logarithmic norm of the residual versus computational cost.}}\label{fig:softmax}
\end{figure}

\begin{figure}[h]
\centering
\begin{subfigure}[b]{0.48\textwidth}
    \includegraphics[width = \textwidth]{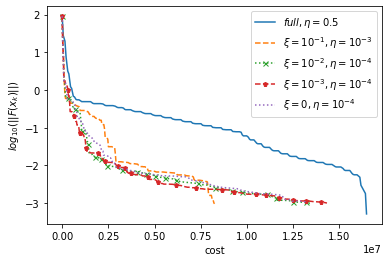}
\end{subfigure}
\begin{subfigure}[b]{0.48\textwidth}
    \includegraphics[width = \textwidth]{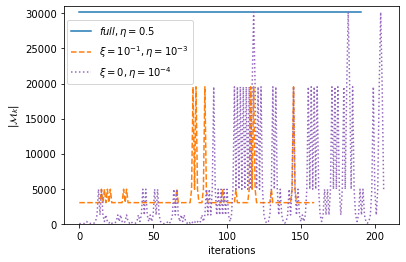}
\end{subfigure}

    \caption{\tblue{Algorithm \metsist and binary classification for the dataset  \emph{a9a}. Exact Jacobians ($\xi=1$) and sparsified Jacobians   with varying  $\xi$ and $\eta$.  Median run in terms of cost: logarithmic norm of the residual versus the computational cost (left) sample size versus the iterations (right).} }
   \label{fig:soft_best}
\end{figure}

On the $x$-axis of Figures \ref{fig:softmax} and \ref{fig:soft_best}   we plot the computational cost, while on the $y$-axis we plot $\|F(x_k)\|$ in logarithmic scale.
Each picture in Figure \ref{fig:softmax} refers to the  median run in terms of computational cost corresponding to a  specific value of $\xi$ and  varying forcing terms. Figures (b), (c) and (d)  display  that, for varying values of $\eta$, the convergence behaviour of the stochastic method is analogous to that of the Gauss-Newton method with full sample shown in Figure (a); on the other hand the stochastic algorithm is computationally more convenient than using the exact Jacobian. Moreover, we note that  the performance of the method is poor for large values of the forcing terms and improves as $\eta$ reduces but then deteriorates again after a certain point. This latter phenomenon is known as {\em oversolving} and indicates that a very accurate solution of the linear systems is pointless \cite{DES}.

In general we know that  as the forcing terms decrease,  the number of outer iterations  decreases as well, while the number of inner iterations for  the solution of  the linear system increases.
Therefore  the most effective choice of the forcing term  depends on the trade-off between the number of inner and outer iterations.
Comparing the subplots in Figure \ref{fig:softmax} we see that
$\eta = 0.5$ is the optimal choice for the case $\xi=1$, i.e., for the   inexact Gauss-Newton method with exact Jacobian, $\eta = 10^{-3}$ is the optimal choice for the case $\xi=0.1$, while for all remaining choices of $\xi$ the forcing term that yields the best result is $\eta = 10^{-4}.$

To summarize, in Figure \ref{fig:soft_best} we compare Algorithm \metsist and the line-search inexact Gauss-Newton method with exact Jacobian. For each value of the scalar $\xi\in\{1, 10^{-1}, 10^{-2}, 10^{-3}, 0\}$, we show the best result in terms of cost obtained varying $\eta$. Figure \ref{fig:soft_best} shows, in particular, that for any value of $\xi<1$, Algorithm \metsist outperforms the  Gauss-Newton method with exact Jacobian ($\xi=1$), and that enforcing a minimum sample size (that is, setting $\xi>0$) can be beneficial. 
In Figure (b), for each run of Figure (a), we plot the value of the sample size $|\mathcal{M}_k|$ along the iterations for the case $\xi=0$, $\xi=10^{-1}$ and for the full Jacobian. Note  that the number of Gauss-Newton iterations performed with $\xi=10^{-1}$ is the lowest and the size of sample size is small for most of the iterations performed. As a result, using $\xi=10^{-1}$ and $\eta=10^{-3}$ provides the best performance. In the case $\xi=0$, the sample sizes is small at many iterations but the number of iterations   is significantly higher than  in the case $\xi=10^{-1}$. Thus,  the per-iteration computational saving that derives from small sample is not sufficient  to balance the increase in the number of iterations, and the run is overall more expensive than  the run with $\xi=10^{-1}$.

\section{Conclusions}
We presented stochastic line-search inexact Gauss-Newton methods for nonlinear least-squares problems and nonlinear systems of equations and analyzed their theoretical properties. Preliminary numerical results indicate that our algorithms are competitive \minor{against} the methods employing exact derivatives.
This work suggests further developments: the generalization of  our algorithms to the case where the residual functions are not evaluated exactly, further investigation of practical rules for fixing the size of the sample, strategies alternative to sampling, such as  sketching techniques, for building the models.


\begin{thebibliography}{99}
\bibitem{bg}  S. Bellavia, G. Gurioli,  
{\em Stochastic analysis of an adaptive cubic regularization method under inexact gradient evaluations and dynamic Hessian accuracy}, Optimization,  71, pp. 227--261, 2022.

\bibitem{BGM} S. Bellavia, G. Gurioli, B. Morini,  {\em Adaptive cubic regularization methods with dynamic inexact Hessian information and applications to finite-sum minimization}, { IMA Journal of Numerical Analysis}, 41, pp. 764--799, 2021.
\bibitem{bkk} S. Bellavia, N. Krejic, N. Krklec Jerinkic, {\em Subsampled Inexact Newton methods for minimizing large sums of convex functions},
{IMA Journal of Numerical Analysis}, 40, pp. 2309--2341, 2020.
%
\bibitem{bkm}  S. Bellavia, N. Kreji\'c, B. Morini, {\em Inexact restoration with subsampled trust-region
methods for finite-sum minimization}, {Computational Optimization and Applications},  {76}, pp. 701--736, 2020. 
%
\bibitem{sirtr_coap}
S. Bellavia, N. Krejić, B. Morini, S. Rebegoldi, {\em A stochastic first-order trust-region method with inexact restoration for finite-sum minimization}, Computational Optimization and Applications, 84, pp. 53–84 2023.

\bibitem{ferrara} 
S. Bellavia,  E. Fabrizi, B. Morini, {\em  Linesearch Newton-CG methods for convex optimization with noise}, Annali dell'Università di Ferrara,  68, pp. 483--504, 2022.

\bibitem {bbn_2017}  {A.S. Berahas, R. Bollapragada, J. Nocedal}, {\em An Investigation of Newton-Sketch and Subsampled Newton Methods}, {Optimization Methods and Software}, 
35, pp.  661--680, 2020. 
%
\bibitem {Berahas}  {A.S. Berahas, L.  Cao, K. Scheinberg}, {\em Global convergence rate analysis of a generic line search algorithm with noise}, {SIAM Journal on Optimization}, 31,  2021.
%
\bibitem{bbn}  R. Bollapragada, R. Byrd, J. Nocedal,  {\em Exact and Inexact Subsampled Newton Methods for
Optimization},  IMA Journal Numerical Analysis, 39, pp. 545--578, 2019.
%
\bibitem{brown} {P.N. Brown and Y. Saad}, Convergence Theory of Nonlinear Newton-Krylov Algorithms {\em SIAM Journal on Optimization, 1994}.
%
\bibitem {cartis}  {C. Cartis, K. Scheinberg}, {\em Global convergence rate analysis of unconstrained optimization methods based on probabilistic model},   Mathematical Programming, 169, pp. 337--375, 2017.
%
\bibitem{STORM1} R. Chen, M. Menickelly, K. Scheinberg,
{\em Stochastic optimization using a trust-region method and random models}, {Mathematical Programming}, {169}, pp. 447-487, 2018.
%
\bibitem{minresqlp} S.C.T. Choi, M.A. Saunders, 
 {\em Algorithm 937: MINRES-QLP for symmetric and Hermitian linear equations and least-squares problems},
ACM Transactions on Mathematical Software, 40, pp. 1–-12, 2014.
%
 \bibitem{DES}  R.S. Dembo, S.C. Eisenstat,  T. Steinhaug, {\em Inexact Newton method}, {SIAM Journal on  Numerical Analysis 19,  
pp. 400--409, 1982.}

\bibitem{LSMR}
D. C.-L. Fong, M. A. Saunders, {\em LSMR: An iterative algorithm for sparse least-squares problems}, SIAM Journal Scientific Computing, 33(5), pp. 2950-2971, 2011.



\bibitem{GR} R.~M. Gower, P. Richtarik, {\em Randomized iterative methods for linear systems}, SIAM Journal on Numerical Analysis, 36(4), pp. 1660-1690, 2015.
\bibitem{HCZ} F. Hegarty, P.\'O Cath\'ain, Y. Zhao, {\em Sparsification of Matrices and Compressed Sensing}, Irish Mathematical Society Bulletin, 81, pp. 5–22, 2018.

\bibitem{HJ} R.A. Horn, C.R. Johnson, {\em Matrix Analysis}, Cambridge University Press, 1985.
\REV{\bibitem{JSX2024}
B. Jin, K. Scheinberg, M. Xie, {\em 
Sample Complexity Analysis for Adaptive Optimization Algorithms with Stochastic Oracles}, Mathematical Programming, 2024.}



\bibitem{integr1} C.T. Kelley, J. I. Northrup, {\em A Pointwise Quasi-Newton Method for Integral Equations}, SIAM Journal on Numerical Analysis, 25, pp. 1138–-1155, 1988.



\bibitem{uci} M. Kelly,  R. Longjohn, K.  Nottingham,
The UCI Machine Learning Repository, https://archive.ics.uci.edu.
%
\bibitem{LN} H. Liu, Q. Ni, {\em Incomplete Jacobian Newton method for nonlinear equations}, An International Journal of Computers \& Mathematics with Applications, 56, pp. 218--227, 2008.
%
\bibitem{LR} Y. Liu, F. Roosta, {\em Convergence of Newton-MR under inexact hessian information}, 
 SIAM Journal on Optimization, 31, pp. 59--90, 2021.
%
\bibitem{acta} P.G. Martinsson, J. A. Tropp, {\em Randomized numerical linear algebra: Foundations and algorithms}, 
Acta Numerica, 29, pp. 403--572, 2020.

%
\bibitem{invex} S.K. Mishra, G. Giorgi, {\em Invexity and optimization}, Vol. 88, Springer Science \& Business Media, 2008.
%

\bibitem{integr2} J.J. Mor\'e, M.Y. Cosnard, {\em Numerical solution of nonlinear equations}, ACM Transactions on  Mathematical  Software,  5, pp. 64--85, 1979.
%
\bibitem{PS}   C. Paquette, K. Scheinberg, {\em A Stochastic Line Search Method with Expected Complexity Analysis},
{SIAM Journal of Optimization}, 30, pp. 349--376, 2020.
%
\bibitem{lsqr} C.C. Paige and M.A. Saunders, {\em LSQR: An algorithm for sparse linear equations and sparse least squares}, ACM Transactions on  Mathematical  Software 8, pp. 
43--71, 1982.
%
\bibitem{RY} L. Reichel, Q. Ye,  {\em Breakdown-free GMRES for singular systems}, SIAM Journal on Matrix Analysis 
and Applications,  26,   pp. 1001--1021, 2005.

\bibitem{Newt-MR1} F. Roosta, Y. Liu, P. Xu, M.W. Mahoney,
{\em Newton-MR: Inexact Newton Method with minimum residual sub-problem solver},  EURO Journal on Computational Optimization, 10, 2022, 100035.


\bibitem{mit1} F. Roosta-Khorasani, M.W. Mahoney, {\em 
Sub-Sampled Newton Methods}, Mathematical  Programming,  174,  pp. 293--326, 2019.
%


\bibitem{ss} Y. Saad, M.H. Schultz, {\em GMRES: a generalized
minimal residual method for solving nonsymmetric linear systems},
SIAM Journal Sci. Stat. Comput., 6 (1985), pp. 856--869.

\bibitem{tropp} J.A. Tropp, {\em  An Introduction to Matrix Concentration Inequalities},
Foundations and Trends in Machine Learning,  8, pp. 1--230, 2015.
%
\bibitem{xu2} P. Xu, F. Roosta-Khorasani, M.W. Mahoney, {\em Newton-Type Methods for Non-Convex Optimization Under Inexact Hessian Information}, Mathematical Programming, 184, pp. 35--70, 2020.


\bibitem{sketched} R. Yuan, A. Lazaric,  R.~M. Gower, {\em  Sketched Newton--Raphson}, SIAM Journal on Optimization, 32, 2022.

\bibitem{wwz2023}
J. Wang, X. Wang, L. Zhang, {\em  Stochastic Regularized Newton Methods for Nonlinear Equations},  Journal of Scientific Computing, 94, article number 51, 2023.

\bibitem{WCK} J. Willert, X. Chen, C~T. Kelley, {\em Newton's method for Monte Carlo-based residuals}, SIAM Journal of Numerical
Analysis, 53, pp. 1738--1757, 2015.
\end{thebibliography}
\end{document}